\documentclass[11pt]{amsart}
\usepackage{amsmath,graphicx}

\newtheorem{thm}{Theorem}[section]
\newtheorem{lem}[thm]{Lemma}%[section]
\newtheorem{cor}[thm]{Corollary}%[section]
\theoremstyle{definition}

\theoremstyle{plain}

\newcommand{\eps}{\varepsilon}
\newcommand{\Z}{{\mathbb{Z}}}

\newcommand{\R}{{\mathbb{R}}}

\newcommand{\C}{{\mathbb{C}}}
\newcommand{\FF}{{\mathcal{F}}}
\newcommand{\T}{{\mathbb{T}}}
\newcommand{\q}{q^\prime}
\newcommand{\N}{{\mathbb{N}}}
\newcommand{\RR}{{\mathcal{R}}}

\newcommand{\X}{\bar{x}}

\newcommand{\EE}{{\mathcal{E}}}

\newcommand{\II}{{\mathcal{I}}}

\newcommand{\JJ}{{\mathcal{J}}}

\newcommand{\ar}{\operatorname{area}}

\newcommand{\dist}{\operatorname{dist}}

\numberwithin{equation}{section}

\begin{document}

\title[Length of a trajectory in a certain billiard in a flat
two-torus]{The average length of a trajectory in a certain
billiard in a flat two-torus}

\author[Boca, Gologan and Zaharescu]{F.P. Boca, R.N. Gologan and A. Zaharescu}

\address{FPB and AZ: Department of Mathematics, University of Illinois at
Urbana-Champaign, 1409 W. Green Street, Urbana, IL 61801, USA}

\address{FPB, RN and AZ: Institute of Mathematics ``Simion Stoilow" of the
Romanian Academy, P.O. Box 1-764, RO-014700 Bucharest, Romania}

\address{E-mail: fboca@math.uiuc.edu; Radu.Gologan@imar.ro; zaharesc@math.uiuc.edu}

\thanks{Research partially supported by ANSTI grant C6189/2000}
\keywords{Periodic Lorentz gas; average first exit time}
\subjclass{11B57; 11P21; 37D50; 58F25; 82C40}

%\date{December 18, 2003}

\begin{abstract}
We remove a small disc of radius $\varepsilon >0$ from the flat
torus ${\mathbb{T}}^2$ and consider a point-like particle that
starts moving from the center of the disk with linear trajectory
under angle $\omega$. Let $\tilde{\tau}_\varepsilon (\omega)$
denote the first exit time of the particle. For any interval
$I\subseteq [0,2\pi)$, any $r>0$, and any $\delta >0$, we estimate
the moments of $\tilde{\tau}_\varepsilon$ on $I$ and prove the
asymptotic formula
\begin{equation*}
\int_I \tilde{\tau}^r_\varepsilon (\omega)\, d\omega = c_r \vert
I\vert \varepsilon^{-r} +O_\delta
(\varepsilon^{-r+\frac{1}{8}-\delta}) \qquad \mbox{\rm as $\
\varepsilon \rightarrow 0^+$},
\end{equation*}
where $c_r$ is the constant
\begin{equation*}
\frac{12}{\pi^2} \int\limits_0^{1/2} \left( x(x^{r-1}+(1-x)^{r-1})
+\frac{1-(1-x)^r}{rx(1-x)} - \frac{1-(1-x)^{r+1}}{(r+1)x(1-x)}
\right) dx.
\end{equation*}
A similar estimate is obtained for the moments of the number of
reflections in the side cushions when ${\mathbb{T}}^2$ is
identified with $[0,1)^2$.
\end{abstract}

\subjclass[2000]{Primary: 11P21 Secondary: 11B57; 37A35; 37A60;
82C05; 82C40}

\maketitle

\section{Introduction and main results}
For each $0<\eps<\frac{1}{2}$ we consider the region
\begin{equation*}
Z_\eps =\{ z\in \R^2 \, ;\, \dist (z,\Z^2) > \eps \}
\end{equation*}
and the {\sl first exit time} (also called {\sl free path length}
by some authors)
\begin{equation*}
\tau_\eps (z,\omega)=\inf \{ \tau >0 \, ;\, z+\tau \omega \in
\partial Z_\eps \},\quad
z\in Z_\eps,\ \omega \in \T ,
\end{equation*}
of a point-like particle which starts moving from the point $z$
with linear trajectory, velocity $\omega$, and constant speed
equal to $1$. This is the model of the periodic two-dimensional
Lorentz gas, intensively studied during the last decades (see
\cite{Ble}, \cite{BGW}, \cite{BS2}, \cite{BS1}, \cite{BSC1},
\cite{BSC2}, \cite{CG}, \cite{Ch1}, \cite{Ch2}, \cite{ChT},
\cite{Dah}, \cite{DDG1}, \cite{DDG2}, \cite{FOK}, \cite{Gal},
\cite{Ku}, \cite{Sin}, \cite{Sin90} for a non-exhaustive list of
references). The phase space of the system consists in the range
of the initial position and velocity and is one of the spaces
$Y_\eps \times \T$ with the normalized Lebesgue measure, or
$\Sigma_\eps^+=\{ (x,y)\in
\partial Y_\eps \times \T\, ;\, \omega \cdot n_x>0\}$ with the
normalized Liouville measure.

\begin{figure}[ht]
\begin{center}
\includegraphics*[scale=0.7, bb=0 0 200 180]{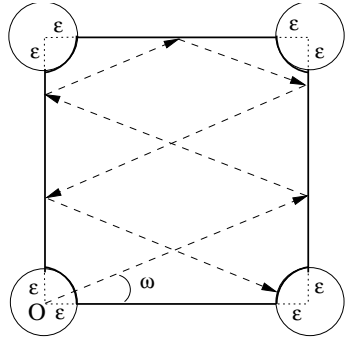}
\end{center}
\caption{The trajectory of the billiard}\label{Figure1}
\end{figure}

Equivalently, one can consider the billiard table $Y_\eps=Z_\eps
/\Z^2$ obtained by removing pockets of the form of quarters of a
circle of radius $\eps$ from the corners. The reflections in the
side cushions are specular and the motion ends when the point-like
particle reaches one of the pockets at the corners. In this
setting $\tau_\eps (z,\omega)$ coincides with the exit time from
the table (see Figure \ref{Figure1}).

This paper considers the situation where the trajectory starts at
the origin $O=(0,0)$. In this case the phase space only consists
in the range of the initial velocity of the particle.  It is given
by the one-dimensional torus $\T$ and can be reduced, for obvious
symmetry reasons, to the interval $\big[ 0,\frac{\pi}{4}\big]$.
From the point of view of Diophantine approximation this
corresponds to a homogeneous problem. We shall be concerned with
estimating the moments of the first exit time $\tilde{\tau}_\eps
(\omega)=\tau_\eps (O,\omega)$ as $\eps \rightarrow 0^+$ when the
phase space is the range $[0,\frac{\pi}{4}]$ of the velocity
$\omega$. This question was raised by Ya.~G. Sinai in a seminar at
the Moskow University in 1981. We answer the question by supplying
asymptotic formulas with explicit main term and error for all the
moments of $\tilde{\tau}_\eps$ in short intervals as follows:

\begin{thm}\label{T1.1}
For any interval $I\subseteq [ 0,\frac{\pi}{4}]$ and any
$r,\delta>0$, one has
\begin{equation*}
\eps^r \int\limits_I \tilde{\tau}_\eps^r (\omega)\, d\omega = c_r
\vert I\vert +\begin{cases} O_{r,\delta}
(\eps^{\frac{1}{8}-\delta}) &
\mbox{\rm if $r\neq 2$} \\
O_{r,\delta} (\eps^{\frac{1}{4}-\delta} ) & \mbox{\rm if $r=2$}
\end{cases} \qquad \mbox{ as $\ \eps \rightarrow 0^+$,}
\end{equation*}
where
\begin{equation*}
c_r =\frac{12}{\pi^2}\int\limits_0^{1/2} \left(
x(x^{r-1}+(1-x)^{r-1}) +\frac{1-(1-x)^r}{rx(1-x)} -
\frac{1-(1-x)^{r+1}}{(r+1)x(1-x)} \right) dx.
\end{equation*}
\end{thm}

The mean free path length is in this case
\begin{equation*}
\frac{4}{\pi} \int\limits_0^{\pi/4} \tilde{\tau}_\eps (\omega)\,
d\omega \sim \frac{c_1}{\eps} =\frac{12}{\pi^2} \cdot \frac{\ln
2}{2\eps} \approx \frac{0.421383}{\eps} .
\end{equation*}
Note also that
\begin{equation*}
\lim\limits_{r\rightarrow 0^+} c_r = -\frac{12}{\pi^2}
\int\limits_0^{1/2} \frac{\ln (1-x)}{x(1-x)}\, dx =1.
\end{equation*}

To prove Theorem \ref{T1.1} we first replace the circular
scatterers by cross-like scatterers $[m-\eps,m+\eps]\times \{n\}
\cup \{ m\} \times [n-\eps,n+\eps]$, $m,n\in \Z^2 \setminus
\{(0,0)\}$.\footnote{Actually it is not hard to see that for
$\omega \in [0,\frac{\pi}{4}]$ the result for cross-like
scatterers is asymptotically the same as when using vertical
scatterers $\{ m\} \times [n-\eps,n+\eps]$.} We denote by $l_\eps
(\omega)$ the free path length in this situation, and first prove

\begin{thm}\label{T1.2}
For any interval $I\subseteq [ 0,\frac{\pi}{4}]$ and any
$r,\alpha,\delta>0$, one has
\begin{equation*}
\eps^r \int\limits_I l_\eps^r (\omega)\, d\omega = c_r
\int\limits_I \frac{dx}{\cos^r x} +\begin{cases} O_{r,\delta}
(\eps^{\frac{1}{2}-2\alpha-\delta}+ \vert I\vert \eps^{\alpha} ) &
\mbox{\rm if $r\neq 2$} \\
O_\delta (\eps^{\frac{1}{2}-\delta}) & \mbox{\rm if $r=2$}
\end{cases} \qquad \mbox{as $\eps \rightarrow 0^+$.}
\end{equation*}
\end{thm}

We consider the probability measures $\tilde{\mu}^I_\eps$ and
$\mu^I_\eps$ on $[0,\infty)$, defined by
\begin{equation*}
\tilde{\mu}^I_\eps (f) =\frac{1}{\vert I\vert} \int\limits_I
f\big( \eps \tilde{\tau}_\eps (\omega) \big) \, d\omega ,\quad
\mu^I_\eps (f) =\frac{1}{\vert I\vert} \int\limits_I f\big( \eps
l_\eps (\omega)\big) \, d\omega ,\qquad f\in C_c ([0,\infty)).
\end{equation*}
Their supports are all contained in $[0,\sqrt{2}]$ as a result of
Lemma \ref{L3.1}. Moreover, we infer from Theorems \ref{T1.1} and
\ref{T1.2} that their moments of order $n\in \N^*$ are of the
form\footnote{We denote $\N=\{ 0,1,2,\dots\}$ and $\N^*=\{
1,2,3,\dots\}$.}
\begin{equation*}
\begin{split}
\tilde{\mu}^I_\eps (X^n) & =\frac{\eps^n}{\vert I\vert}
\int\limits_I \tilde{\tau}^n_\eps (\omega)\, d\omega =c_n
+\frac{1}{\vert I\vert} \, O_{n,\delta}
(\eps^{\frac{1}{8}-\delta});\\
\mu^I_\eps (X^n) & = \frac{\eps^n}{\vert I\vert} \int\limits_I
l_\eps^n (\omega)\, d\omega=\frac{c_n}{\vert I\vert} \int\limits_I
\frac{dx}{\cos^n x} +\frac{1}{\vert I\vert}\, O_{n,\delta} (
\eps^{\frac{1}{6}-\delta} ).
\end{split}
\end{equation*}
These asymptotic formulas show in particular that
$\tilde{\mu}^I_\eps (X^n)$ and $\mu^I_\eps (X^n)$ converge to the
main terms as $\varepsilon \rightarrow 0^+$. The Banach-Alaoglu
and Stone-Weierstrass theorems now lead to

\begin{cor}\label{C1.3}
There exist probability measures $\tilde{\mu}$ and $\tilde{\mu}^I$
on $[0,\sqrt{2}]$ such that
\begin{equation*}
\tilde{\mu}^I_\eps \rightarrow \tilde{\mu} \quad \mbox{and} \quad
\mu^I_\eps \rightarrow \mu^I \quad \mbox{weakly as $\ \eps
\rightarrow 0^+$}.
\end{equation*}
Moreover, the moments of $\tilde{\mu}$ and $\mu^I$ are
\begin{equation*}
\int\limits_0^\infty t^n \, d\tilde{\mu}(t)=c_n
\end{equation*}
and respectively
\begin{equation*}
\int\limits_0^\infty t^n \, d\mu^I (t)= \frac{c_n}{\vert I\vert}
\int\limits_I \frac{dx}{\cos^n x} \, .
\end{equation*}
\end{cor}

Besides, we estimate the average of the number of reflections
$\tilde{R}_\eps (\omega)$ in the side cushions of the billiard
table in the case of circular scatterers and prove
\begin{thm}\label{T1.4}
For any interval $I\subseteq [0,\frac{\pi}{4}]$ and any
$r,\delta>0$, one has
\begin{equation*}
\eps^r \int\limits_I \tilde{R}_\eps^r (\omega)\, d\omega =c_r
\int\limits_I (\sin x+\cos x)^r \, dx+O_{r,\delta}
(\eps^{\frac{1}{8}-\delta}) \quad \mbox{ as $\ \eps \rightarrow
0^+$}.
\end{equation*}
\end{thm}

Again, we first consider the case of cross-like (or vertical)
scatterers, let $R_\eps (\omega)$ denote the number of reflections
in the side cushions of the billiard table in this case, and prove

\begin{thm}\label{T1.5}
For any interval $I\subseteq [ 0,\frac{\pi}{4}]$ and any
$r,\alpha,\delta>0$, one has
\begin{equation*}
\eps^r \int\limits_I R_\eps^r (\omega)\, d\omega = c_r
\int\limits_I (1+\tan x)^r \, dx + O_{r,\delta}
(\eps^{\frac{1}{2}-2\alpha-\delta} +\vert I\vert \eps^{\alpha})
\qquad \mbox{as $\eps \rightarrow 0^+$.}
\end{equation*}
\end{thm}

We may also consider the probability measures $\tilde{\nu}^I_\eps$
and $\nu^I_\eps$ on $[0,\infty)$ associated with the random
variables $\eps \tilde{R}_\eps$ and $\eps R_\eps$, and defined by
\begin{equation*}
\tilde{\nu}^I_\eps (f) =\frac{1}{\vert I\vert} \int\limits_I
f\big( \eps \tilde{R}_\eps (\omega) \big) \, d\omega ,\quad
\nu^I_\eps (f) =\frac{1}{\vert I\vert} \int\limits_I f\big( \eps
R_\eps (\omega)\big) \, d\omega ,\qquad f\in C_c ([0,\infty)).
\end{equation*}
From Theorems \ref{T1.4} and \ref{T1.5} we derive

\begin{cor}\label{C1.6}
There exist probability measures $\tilde{\nu}^I$ and
$\tilde{\nu}^I$ on $[0,\sqrt{2}]$ such that
\begin{equation*}
\tilde{\nu}^I_\eps \rightarrow \tilde{\nu}^I \quad \mbox{and}
\quad \nu^I_\eps \rightarrow \nu^I \quad \mbox{as $\ \eps
\rightarrow 0^+$}.
\end{equation*}
Moreover, the moments of $\tilde{\nu}^I$ and $\nu^I$ are
\begin{equation*}
\int\limits_0^\infty t^n \, d\tilde{\nu}^I (t) =\frac{c_n}{\vert
I\vert} \int\limits_I (\sin x+\cos x)^n \, dx,
\end{equation*}
and respectively
\begin{equation*}
\int\limits_0^\infty t^n \, d\nu^I (t)= \frac{c_n}{\vert I\vert}
\int\limits_I (1+\tan x)^n \, dx.
\end{equation*}
\end{cor}

In the case $I\subseteq [\frac{\pi}{4},\frac{\pi}{2}]$ one gets
formulas similar to the ones in Theorems \ref{T1.1}, \ref{T1.2},
\ref{T1.4} and \ref{T1.5} after performing a symmetry with respect
to a diagonal of the square, i.e. replacing $(\alpha,\beta)$ by
$(\frac{\pi}{2}-\beta,\frac{\pi}{2}-\alpha)$.

The proofs make use of techniques employed in the study of the
spacings between Farey fractions, pioneered in \cite{Hall},
\cite{Hall2}, \cite{HT}, and furthered recently in \cite{BCZ1},
\cite{BCZ2}, \cite{ABCZ}, \cite{HZ} where estimates for
Kloosterman sums are being used. The first step consists in
proving Theorems \ref{T1.2} and \ref{T1.5}, which refer to the
case of cross-like or vertical scatterers. In this case one can
directly take advantage of the fact that the intervals $I_{a/q}
=[\frac{a-\eps}{q},\frac{a+\eps}{q}]$, with $\frac{a}{q}$ Farey
fraction of order $Q=[\frac{1}{\eps}]$, provide a covering of
$[0,1]$ such that two intervals $I_{a/q}$ and
$I_{a^\prime/q^\prime}$ overlap if and only if $\frac{a}{q}$ and
$\frac{a^\prime}{q^\prime}$ are consecutive Farey fractions of
order $Q$.

Finally, the case of circular scatterers is settled by
partitioning the range $I$ into $[\eps^{-\theta}]$ intervals of
equal size for a convenient value of the exponent $\theta$, and
replacing the small circles of radius $\eps$ first by vertical
scatterers of type $\{ m\} \times [n-\eps_- (m,n),n+\eps_+
(m,n)]$, and finally by scatterers of type $\{ m\} \times
[n-\tilde{\eps},n+\tilde{\eps}]$ for appropriate choices of
$\eps_\pm (m,n)$ and $\tilde{\eps}$.

It should be possible in theory to compute the densities of the
limit measures from their moments using either the Cauchy
transform or the inverse Mellin transform. An attempt of this kind
does not seem to easily lead however to a tractable formula for
these densities. The convergence of the measures
$\tilde{\mu}^I_\eps$ and $\tilde{\nu}^I_\eps$ was proved in a
different way and the limit measures were explicitly computed in
\cite{BGZ}.

Techniques using Farey fractions and Kloosterman sums were
recently used in \cite{BZ} to establish the existence, and compute
the distribution, of the free path length for the periodic
two-dimensional Lorentz gas in the small-scatterer limit in the
case where the trajectory does not necessary start from the
origin, and one averages over both initial position and initial
velocity.

This is the final version of the paper with the same title,
circulated as preprint math.NT/0110208.

\bigskip

\section{Farey fractions and Kloosterman sums}

For each integer $Q\geq 1$, let $\FF_Q$ denote the set of Farey
fractions of order $Q$, i.e. irreducible fractions in the interval
$(0,1]$ with denominator $\leq Q$. The number of Farey fractions
of order $Q$ in an interval $J\subseteq [0,1]$ can be expressed as
\begin{equation*}
\# (J\cap \FF_Q) =\frac{Q^2\vert J\vert}{2\zeta (2)} +O(Q\ln Q) .
\end{equation*}

Recall that if $\frac{a}{q}<\frac{a^\prime}{\q}$ are two
consecutive elements in $\FF_Q$, then
\begin{equation*}
a^\prime q-a\q =1 \qquad \mbox{\rm and} \qquad q+\q >Q.
\end{equation*}
Conversely, if $q,\q \in \{1,\dots ,Q\}$ and $q+\q>Q$, then there
are $a\in \{ 1,\dots,q-1\}$ and $a^\prime \in \{ 1,\dots, \q-1\}$
such that $\frac{a}{q} < \frac{a^\prime}{\q}$ are consecutive
elements in $\FF_Q$. Proofs of these well-known properties of
Farey fractions can be found for instance in \cite{HW},
\cite{Hall}, \cite{LV}.

Throughout the paper we shall denote by $\FF_Q^<$, and
respectively by $\FF_Q^>$, the set of pairs
$(\frac{a}{q},\frac{a^\prime}{q^\prime})$ of consecutive elements
in $\FF_Q$ with $q<q^\prime$, and respectively with $q>q^\prime$.
We also set
\begin{equation*}
\begin{split}
& \Z^2_{\mathrm{pr}} =\{ (a,b) \in \Z^2 \, ; \, \gcd (a,b)=1\}; \\
&\sideset{^J}{} \sum\limits_{a/q} =\sum_{\substack{(
a/q,a^\prime/\q) \in \FF_Q^< \\ a/q \in J}} \qquad \mbox{\rm and}
\qquad \sideset{}{^J} \sum\limits_{a/q} =\sum_{\substack{(
a/q,a^\prime/\q ) \in \FF_Q^> \\ a/q \in J}} ; \\ & \Delta_Q =\{
(x,y)\in \Z^2_{\mathrm{pr}} \, ;\, 0<x,y\leq Q ,\ x+y>Q \} ; \\ &
\RR_{m,n} =[m,m+1] \times [n,n+1], \qquad m,n\in \R .
\end{split}
\end{equation*}

For each region $\RR$ in $\R^2$ and each $C^1$ function $f:\RR
\rightarrow \C$, we denote
\begin{equation*}
\| f\|_{\infty ,\RR} =\sup\limits_{(x,y)\in \RR} \vert f(x,y)\vert
,\quad \| Df\|_{\infty ,\RR} =\sup\limits_{(x,y)\in \RR} \bigg(
\bigg| \frac{\partial f}{\partial x} \, (x,y)\bigg| + \bigg|
\frac{\partial f}{\partial y} \, (x,y) \bigg| \bigg) .
\end{equation*}

The notation $f\ll g$ means the same thing as $f=O(g)$; that is,
there exists an absolute constant $c>0$ such that $\vert f\vert
\leq cg$ for all values of the variable under consideration. When
the constant depends on a parameter $\delta$, this dependence will
be indicated by writing $f\ll_\delta g$. The notation $f\asymp g$
will mean that $f\ll g$ and $g\ll f$ simultaneously.

We shall be mainly interested in consecutive Farey fractions
$\frac{a}{q} <\frac{a^\prime}{\q}$ in $\FF_Q$ with the property
that, say, $\frac{a}{q}$ belongs to a prescribed interval
$J\subseteq [0,1]$. The equality $a^\prime q-a\q=1$ yields
$a=q-\bar{\q}$, where $\X$ denotes the unique integer in $\{
1,2,\dots ,q-1\}$ for which $x\X=1\hspace{-2pt} \pmod{q}$. Thus
$\frac{a}{q} \in J=[t_1,t_2 ]$ is equivalent to $\bar{\q} \in
J_q^{(1)} :=[ (1-t_2)q,(1-t_1 )q]$. Moreover, $\frac{a^\prime}{\q}
\in J$ is equivalent to $\bar{q} \in J_{\q}^{(2)} :=[t_1 \q ,t_2
\q]$, where this time $\bar{q}$ denotes the multiplicative inverse
of $q\hspace{-2pt} \pmod{\q}$.

An important device employed in \cite{BCZ1}, \cite{BCZ2},
\cite{ABCZ} to estimate sums over primitive lattice points is the
Weil type \cite{We} estimate
\begin{equation}\label{2.1}
\vert S(m,n;q)\vert \ll \tau (q) \gcd (m,n,q)^{\frac{1}{2}}
q^{\frac{1}{2}}
\end{equation}
on complete Kloosterman sums
\begin{equation*}
S(m,n;q)=\sum_{\substack{x\in [1,q] \\ \gcd (x,q)=1}} e\bigg(
\frac{mx+n\bar{x}}{q} \bigg),
\end{equation*}
in the presence of an integer albeit not necessarily prime modulus
$q$, proved in \cite{Hoo} (see also \cite{Est}). The bound from
\eqref{2.1} can be used (see \cite[Lemma\,1.7]{BCZ2}) to prove the
estimate
\begin{equation}\label{2.2}
N_q (\II,\JJ)=\frac{\varphi (q)}{q^2} \, \vert \II \vert \, \vert
\JJ \vert+O_\delta ( q^{\frac{1}{2}+\delta} )
\end{equation}
for the number $N_q (\II,\JJ)$ of pairs of integers $(x,y)\in \II
\times \JJ$ for which $xy=1\hspace{-2pt}\pmod{q}$, whenever $\II$
and $\JJ$ are intervals which contain at most $q$ integers.

We shall use the following slight improvement of Corollary 1 and
Lemma 8 in \cite{BCZ2}. The proof follows literally the reasoning
from Lemmas 2, 3 and 8 in \cite{BCZ2}.

\begin{lem}\label{L2.1}
Let $\Omega \subseteq [1,R] \times [1,R]$ be a convex region and
let $f$ be a $C^1$ function on $\Omega$. Then

\medskip

{\em (i)} $\quad \displaystyle \sum\limits_{(a,b)\in \Omega \cap
\Z^2_{\mathrm{pr}}} \hspace{-2pt} f(a,b) =\frac{1}{\zeta (2)}
\iint\limits_\Omega f(x,y) \, dx\, dy +\EE_{R,\Omega,f},$

where
\begin{equation*}
\EE_{R,\Omega,f} \ll \| f\|_{\infty,\Omega} R \ln R+ \hspace{-2pt}
\sum_{\substack{(m,n)\in \Z^2 \\ \RR_{m,n} \subset
\overline{\Omega}}} \hspace{-2pt} \| Df\|_{\infty, \RR_{m,n}}
\hspace{-1.5pt} \ln R .
\end{equation*}

\medskip

{\em (ii)} For any interval $J\subseteq [0,1]$ one has
\begin{equation*}
\sum_{\substack{(a,b)\in \Omega \cap \Z^2_{\mathrm{pr}} \\
\bar{b} \in J_a}} \hspace{-2pt} f(a,b)=\frac{\vert J\vert}{\zeta
(2)} \iint\limits_\Omega \hspace{-1.5pt} f(x,y)\, dx\, dy
+\FF_{R,\Omega,f,J},
\end{equation*}
where
\begin{equation*}
\FF_{R,\Omega,f,J} \ll_\delta \| f\|_{\infty,\Omega} m_f
R^{\frac{3}{2}+\delta} +\| f\|_{\infty,\Omega} \mbox{\rm length}
(\partial \Omega) \ln R+\hspace{-8pt} \sum_{\substack{(m,n)\in
\Z^2 \\ \RR_{m,n} \in \overline{\Omega}}} \hspace{-5pt} \|
Df\|_{\infty,\RR_{m,n}} \ln R
\end{equation*}
for any $\delta >0$, where $\bar{b}$ denotes\footnote{When writing
$\bar{b}\in J_a$ we implicitly assume that $\gcd (a,b)=1$.} the
multiplicative inverse of $b\hspace{-2pt} \pmod{a}$, $J_a$ is
either $J_a^{(1)}$ or $J_a^{(2)}$, and $m_f$ is an upper bound for
the number of intervals of monotonicity of each of the functions
$y\mapsto f(x,y)$.
\end{lem}

The proof of $(ii)$ relies on \eqref{2.2}. We also note the
following important corollary of \eqref{2.2}, which will be often
employed in this paper and in the subsequent work from \cite{BGZ}
and \cite{BZ}.

\begin{lem}\label{L2.2}
Assume that $q\geq 1$ is an integer, $\II$ and $\JJ$ are intervals
which contain at most $q$ integers, and $f:\II \times \JJ
\rightarrow \R$ is a $C^1$ function. Then for any integer $T>1$
one has
\begin{equation*}
\sum_{\substack{a\in \II,\, b\in \JJ \\
ab=1\hspace{-8pt}\pmod{q}}} \hspace{-3pt} f(a,b)  =\frac{\varphi
(q)}{q^2} \iint\limits_{\II \times \JJ} f(x,y) dx dy
+E_{q,\II,\JJ,f,T},
\end{equation*}
where
\begin{equation*}
E_{q,\II,\JJ,f,T} \ll_\delta T^2 q^{\frac{1}{2}+\delta} \|
f\|_\infty+Tq^{\frac{3}{2}+\delta} \| Df\|_\infty+ \frac{\vert \II
\vert \, \vert \JJ \vert \, \| Df\|_\infty}{T}
\end{equation*}
for all $\delta >0$. Here $\| \cdot \|_\infty$ denotes the
$L^\infty$-norm on $\II \times \JJ$.
\end{lem}

\begin{proof} If $T\geq q$, then the error is larger than the sum to
estimate and there is nothing to prove.

If $T<q$, we partition the intervals $\II$ and $\JJ$ respectively
into $T$ intervals $\II_1,\dots,\II_{\mbox{\tiny $T$}}$ and
$\JJ_1,\dots,\JJ_{\mbox{\tiny $T$}}$ of equal size $\vert \II_i
\vert =\frac{\vert \II\vert}{T}$ and $\vert \JJ_j\vert=\frac{\vert
\JJ\vert}{T}$. The idea is to approximate $f(x,y)$ by a constant
whenever $(x,y)\in \II_i \times \JJ_j$. For, we choose for each
pair of indices $(i,j)$ a point $(x_{ij},y_{ij}) \in \II_i \times
\JJ_j$ for which
\begin{equation}\label{2.3}
\iint\limits_{\II_i \times \JJ_j} f=\vert \II_i \vert \, \vert
\JJ_j \vert f(x_{ij},y_{ij} ).
\end{equation}

For $(x,y)\in \II_i \times \JJ_j$ the mean value theorem gives
\begin{equation}\label{2.4}
\begin{split}
f(x,y) & =f(x_{ij},y_{ij})+ O\big( (\vert \II_i \vert +\vert
\JJ_j\vert)\| Df\|_\infty \big) \\ & =f(x_{ij},y_{ij})+ O\bigg(
\frac{q}{T}\, \| Df\|_\infty \bigg) .
\end{split}
\end{equation}
This gives in turn
\begin{equation}\label{2.5}
\begin{split}
\sum_{\substack{a\in \II ,b\in \JJ \\ ab=1\hspace{-8pt}\pmod{q}}}
f(a,b) & =\sum\limits_{i,j=1}^T \sum_{\substack{(x,y)\in \II_i
\times \JJ_j \\ xy=1\hspace{-8pt} \pmod{q}}} f(x,y) \\ &
=\sum\limits_{i,j=1}^T N_q (\II_i,\JJ_j)\left(
f(x_{ij},y_{ij})+O\Big( \frac{q\| Df\|_\infty}{T}\Big) \right).
\end{split}
\end{equation}

Since each interval $\II_i$ and $\JJ_j$ contains at most $q$
integers, estimate \eqref{2.2} applies to them and gives
\begin{equation}\label{2.6}
N_q (\II_i ,\JJ_j)=\frac{\varphi(q)}{q^2} \, \vert \II_i \vert \,
\vert \JJ_j \vert +O_\delta (q^{\frac{1}{2}+\delta}).
\end{equation}

As a result of \eqref{2.6} and \eqref{2.3}, the main term in
\eqref{2.5} becomes
\begin{equation*}
\begin{split}
& \frac{\varphi(q)}{q^2}\sum\limits_{i,j=1}^T \vert \II_i\vert \,
\vert \JJ_j\vert f(x_{ij},y_{ij})+O_\delta (T^2
q^{\frac{1}{2}+\delta} \| f\|_\infty ) \\
& \qquad =\frac{\varphi(q)}{q^2} \int\limits_{\II \times \JJ} f\
+O_\delta (T^2 q^{\frac{1}{2}+\delta} \| f\|_\infty),
\end{split}
\end{equation*}
while the error term in \eqref{2.5} will be
\begin{equation*}
\ll \frac{q\| Df\|_\infty}{T} \bigg( \frac{\varphi(q)}{q^2}\,
\vert \II\vert\, \vert \JJ \vert +T^2 q^{\frac{1}{2}+\delta}\bigg)
\leq \| Df\|_\infty \bigg( \frac{\vert \II\vert \, \vert
\JJ\vert}{T}+Tq^{\frac{3}{2}+\delta}\bigg).
\end{equation*}
\end{proof}

\bigskip

\section{The second moment of the first exit time for
cross-like scatterers}

Throughout this section we keep $0<\eps <\frac{1}{2}$ fixed, and
take
\begin{equation*}
Q=Q_\eps=\left[ \frac{1}{\eps}\right] = \mbox{\rm the integer part
of}\ \frac{1}{\eps}\, .
\end{equation*}
We also denote
\begin{equation*}
\begin{split}
\Z^{2 \ast} & =\Z^2 \setminus \{ (0,0)\} ,\\
C_\eps & =\{ 0\} \times
[-\eps,\eps] \cup [-\eps,\eps ]\times \{ 0\} ,\\
V_\eps & =\{ 0\} \times [-\eps,\eps], \\
l_\eps (\omega) & =\inf \{ \tau>0 \, ;\, (\tau \cos \omega,\tau
\sin \omega)
\in C_\eps +\Z^{2 \ast} \} \\
t_P & =\mbox{\rm the slope of the line $OP$,} \\
\| (x,y)\| & =\sqrt{x^2+y^2}\, ,  \qquad x,y\in \R .
\end{split}
\end{equation*}
Let
\begin{equation*}
{\mathfrak C}_\eps =C_\eps +\{ (q,a) \, ;\, a/q \in \FF_Q \}
\end{equation*}
denote the translates of $C_\eps$ at all integer points with slope
in $\FF_Q$.

For each point $A(q,a)$ with $\frac{a}{q} \in \FF_Q$ we construct
a vertical segment $NS$ of length $2\eps$ and a horizontal segment
$WE$ of length $2\eps$, both centered at $A$.

Performing symmetries with respect to the integer vertical and
horizontal lines, the problem translates into a covering version
in $\R^2$. It is clear that one can discard the points
$(q^\prime,a^\prime)$ with $\gcd (q^\prime ,a^\prime)=d>1$, which
are already hidden by $\big( \frac{q^\prime}{d},\frac{a^\prime}{d}
\big)$.

The trajectory will now originate at $O=(0,0)$ and end when it
reaches one of the components $(q,a)+C_\eps$ of ${\mathfrak
C}_\eps$, $\frac{a}{q}\in \FF_Q$, as seen in the next elementary
but useful lemma.

\begin{lem}\label{L3.1}
Any ray of direction $\omega \in [0,\frac{\pi}{4}]$ which
originates at $O$ inevitably intersects ${\mathfrak C}_\eps$.
Moreover, if
$\gamma=\frac{a}{q}<\gamma^\prime=\frac{a^\prime}{\q}$ are two
consecutive Farey fractions in $\FF_Q$ and $\tan \omega \in [
\gamma , \gamma^\prime ]$, then the ray of direction $\omega$
intersects either $(q,a)+C_\eps$ or $(q^\prime,a^\prime)+C_\eps$
and does not intersect any other component of ${\mathfrak
C}_\eps$.\footnote{Equivalently, the intervals
$I_{a/q}=[\frac{a-\eps}{q},\frac{a+\eps}{q}]$, $\frac{a}{q} \in
\FF_Q$, cover $[0,1]$ and two such intervals $I_{a/q}$ and
$I_{a^\prime/q^\prime}$ overlap if and only if $\frac{a}{q}$ and
$\frac{a^\prime}{q^\prime}$ are consecutive elements in $\FF_Q$.}
\end{lem}

\begin{proof} We shall utilize the inequalities $q+\q \geq
Q+1>\frac{1}{\eps}\geq Q\geq \max \{ q,\q \}$, getting
\begin{equation*}
t_A =\frac{a}{q} \leq t_{S^\prime}=\frac{a^\prime-\eps}{\q} <t_N
=\frac{a+\eps}{q} \leq t_{A^\prime} =\frac{a^\prime}{\q} \, .
\end{equation*}

\begin{figure}[ht]
\includegraphics*[scale=1, bb=0 0 210 140]{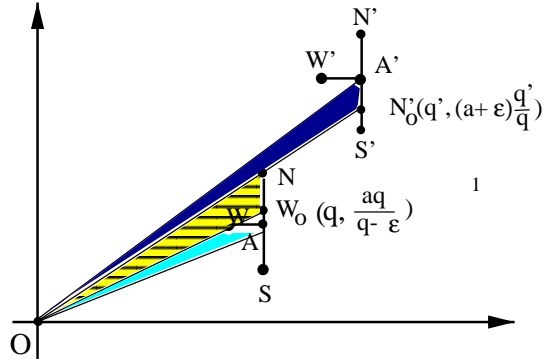}
\caption{The case $q<\q$}\label{Figure2}
\end{figure}

\begin{figure}[ht]
\includegraphics*[scale=0.7, bb=0 0 300 200]{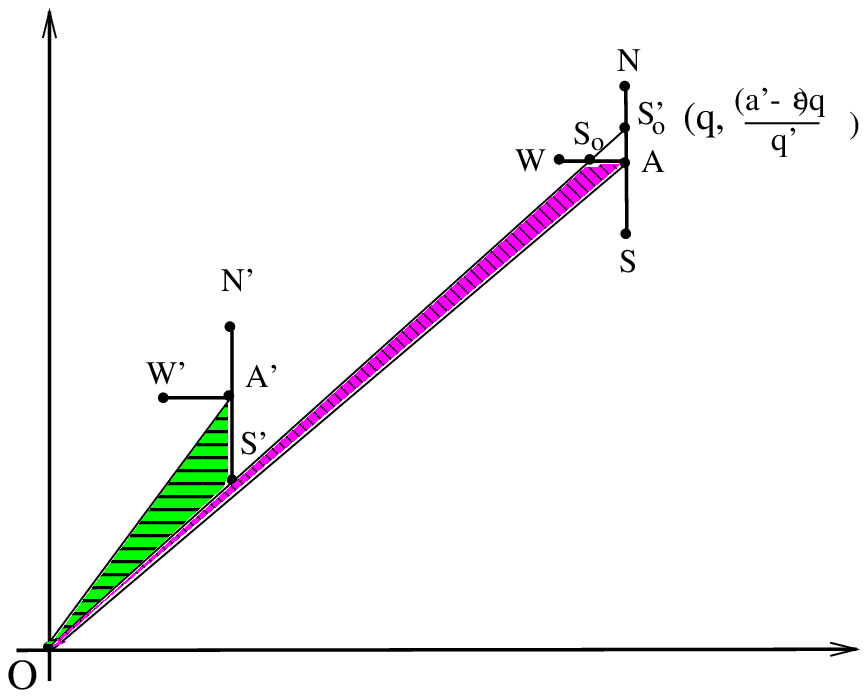}
\includegraphics*[scale=0.75, bb=0 0 230 190]{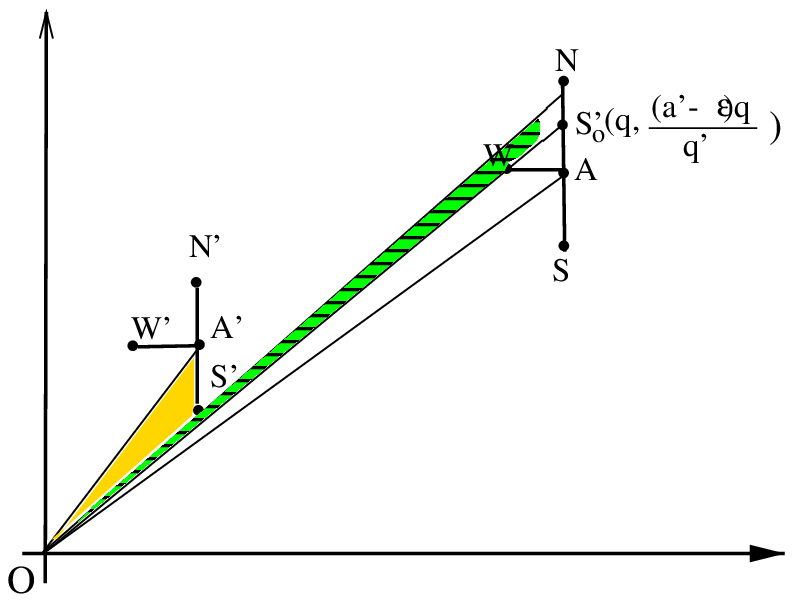}
\caption{The case $\q <q$ and $t_{S^\prime} \leq t_W$,
respectively $\q<q$ and $t_{S^\prime} >t_W$.} \label{Figure3}
\end{figure}

In the case $q<\q$, we set $\{ W_0\}=OW \cap NS$ and $\{N_0^\prime
\}=ON \cap N^\prime S^\prime$ (see Figure \ref{Figure2} and note
that $a<a^\prime$), inferring that
\begin{equation}\label{3.1}
\begin{split}
\int\limits_{\arctan \gamma}^{\arctan \gamma^\prime}
\hspace{-12pt} l_\eps^2 (\omega) \, d\omega & =2\ar (\triangle
OAN)+2\ar (\triangle ON_0^\prime A^\prime) -2\ar (\triangle AW_0
W) \\ & =\eps q+\q \bigg( a^\prime-\frac{(a+\eps)\q}{q}
\bigg)+O(\eps^2) \\ & =\frac{\q-\eps (q^{\prime 2}-q^2)}{q}+O(
\eps^2) .
\end{split}
\end{equation}

In the case $q>\q$ one has $a^\prime<a$. Moreover,
\begin{equation*}
t_A =\gamma \leq \min (t_W,t_{S^\prime})=\min \bigg(
\frac{a}{q-\eps}\, ,\frac{a^\prime-\eps}{\q} \bigg) \leq \max \{
t_W,t_{S^\prime} \} \leq \gamma^\prime =t_{A^\prime} .
\end{equation*}
This shows that any ray of slope $\tan \omega \in
[\gamma,\gamma^\prime]$ intersects either $(q,a)+C_\eps$ or
$(\q,a^\prime)+C_\eps$ and no other component of ${\mathfrak
C}_\eps$ (see Figures \ref{Figure2} and \ref{Figure3}).

Besides, we estimate the average of the second moment of the
length $l_\eps (\omega)$ of the trajectory when $\tan \omega \in
[\gamma,\gamma^\prime]$.

When $t_{S^\prime} \leq t_W$ (i.e.
$a^\prime+q>\eps+\frac{1}{\eps}$), we set $\{ S_0 \}=OS^\prime
\cap AW$, $\{ S_0^\prime \}=OS^\prime \cap NS$, and note (see the
first picture in Figure \ref{Figure3}) that
\begin{equation}\label{3.2}
\begin{split}
\int\limits_{\arctan \gamma}^{\arctan \gamma^\prime} \hspace{-8pt}
l_\eps^2 (\omega) \, d\omega & =2 \ar (\triangle OA^\prime
S^\prime)+ 2\ar (\triangle OAS_0^\prime)-2\ar (\triangle AS_0
S_0^\prime) \\ & =\eps \q+q\bigg( \frac{(a^\prime-\eps)q}{\q}-a
\bigg) +O(\eps^2) \\ & =\frac{q-\eps (q^2 -q^{\prime 2})}{\q}
+O(\eps^2).
\end{split}
\end{equation}

When $t_{S^\prime}>t_W$, we set $\{ W_0 \} =OW \cap NS$, $\{
S_0^\prime \}=OS^\prime \cap NS$, and get (see the second picture
in Figure \ref{Figure3})
\begin{equation}\label{3.3}
\begin{split}
\int\limits_{\arctan \gamma}^{\arctan \gamma^\prime} \hspace{-8pt}
l_\eps^2 (\omega) \, d\omega & =2\ar (\triangle OA^\prime
S^\prime)+ 2\ar (\triangle OAS_0^\prime)-2\ar (\triangle AWW_0 )
\\ & =\frac{q-\eps (q^2 -q^{\prime 2})}{\q}
+O(\eps^2).
\end{split}
\end{equation}
\end{proof}

We consider the region
\begin{equation*}
\Omega_Q=\{ (x,y)\in \R^2 \, ;\, 1\leq x\leq y\leq Q,\ x+y>Q\}
\end{equation*}
and the function
\begin{equation*}
f(x,y)=\frac{y+\varepsilon (x^2-y^2)}{x} =\frac{y(1-\varepsilon
y)}{x}+\varepsilon x, \qquad (x,y)\in \Omega_Q .
\end{equation*}

\begin{figure}[t]
\includegraphics*[scale=0.9, bb=0 0 180 180]{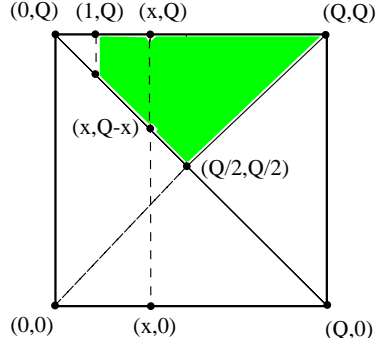}
\caption{The region $\Omega_Q$} \label{Figure4}
\end{figure}

Consider also $I=[\alpha ,\beta ]\subseteq [ 0,\frac{\pi}{4}]$,
take $t_1 =\tan \alpha$, $t_2 =\tan \beta$, and let $J=[t_1 ,t_2 ]
\subseteq [0,1]$. For $(x,y)\in \Omega_Q$ one has $x>Q-y>
\frac{1}{\eps}-1-y$, which gives $1-\eps y<\eps (x+1)\leq 2\eps
x$. It is also seen that $1-\eps y\geq 1-\frac{y}{Q}\geq 0$. As a
result we find that $\| f\|_{\infty,\Omega_Q} \leq 3$. Since
$\varepsilon^2 \# \FF_Q < 1$, formulas \eqref{3.1}, \eqref{3.2},
\eqref{3.3} provide
\begin{equation}\label{3.4}
\int\limits_I l_\eps^2 (\omega) \, d\omega =2\sideset{^J}{}
\sum\limits_{a/q} \ \int\limits_{\arctan \frac{a}{q}}^{\arctan
\frac{a^\prime}{\q}} \hspace{-8pt} l_\eps^2 (\omega) \, d\omega
+O(1) =2\sideset{^J}{} \sum\limits_{a/q} f(q,\q) +O(1).
\end{equation}

To master the latest sum, we aim to apply Lemma \ref{L2.1} to
$\Omega_Q$. With the notation from Section 2, relation \eqref{3.4}
yields
\begin{equation}\label{3.5}
\int\limits_I l_\eps^2 (\omega) \, d\omega =2
\sum_{\substack{(a,b)\in \Omega_Q \\ \bar{b} \in J_a^{(1)}}}
\hspace{-4pt} f(a,b)+O(1).
\end{equation}

We also see that for $(x,y)\in \Omega_Q$ one has $\big|
\frac{\partial f}{\partial x}\big| =\big| \eps-\frac{y(1-\eps
y)}{x^2}\big| \leq \eps +\frac{2\eps y}{x} \leq \frac{3}{x}$ and
$\big| \frac{\partial f}{\partial y}\big| =\frac{\vert 1-2\eps
y\vert}{x} \leq \frac{1}{x}$; hence
\begin{equation}\label{3.6}
\sum\limits_{\substack{(a,b)\in \Z^2 \\
\RR_{a,b}\subset \bar{\Omega}_Q}} \| Df\|_{\infty ,\RR_{a,b}} \ll
\sum\limits_{x=1}^Q \ \sum\limits_{y=\max \{Q-x,x\}}^Q \frac{1}{x}
\ \ll \sum\limits_{x=1}^Q 1=Q.
\end{equation}

Now we can apply Lemma \ref{L2.1} (ii) to the sum from
\eqref{3.5}, and employ \eqref{3.6} and $m_f \leq 2$, to infer
that
\begin{equation}\label{3.7}
\int\limits_I l_\eps^2 (\omega) \, d\omega
=\frac{2(t_2-t_1)}{\zeta(2)} \iint\limits_{\Omega_Q} f(x,y)\, dx\,
dy+O_\delta (Q^{\frac{3}{2}+\delta} ) .
\end{equation}

When $\alpha =0$ and $\beta=\frac{\pi}{4}$, Lemma \ref{L2.1} $(i)$
improves upon the error in \eqref{3.7} to
\begin{equation}\label{3.8}
\int\limits_0^{\pi/4} l_\eps^2 (\omega) \,
d\omega=\frac{2}{\zeta(2)} \iint\limits_{\Omega_Q} f(x,y)\,
dx\,dy+O(Q\ln Q).
\end{equation}

In summary, \eqref{3.7}, \eqref{3.8} and the equality
\begin{equation*}
\iint\limits_{\Omega_Q} f(x,y)\, dx\,dy=\frac{1+2\ln 2}{12} \, Q^2
\end{equation*}
lead to
\begin{thm}\label{T3.2}
$\mathrm{(i)} \quad \displaystyle \eps^2 \int\limits_0^{\pi/4}
l_\eps^2 (\omega) \, d\omega =\frac{1+2\ln 2}{\pi^2} +O( \eps
\vert \ln \eps \vert ) \qquad \mbox{as $\ \eps \rightarrow 0^+$.}$

$\mathrm{(ii)}$ For any $0\leq \alpha<\beta \leq \frac{\pi}{4}$
and $\delta >0$, one has
\begin{equation*}
\eps^2 \int\limits_\alpha^\beta l_\eps^2 (\omega)\, d\omega =
\frac{(1+2\ln 2)(\tan \beta -\tan \alpha)}{\pi^2} +O_\delta
(\eps^{\frac{1}{2}-\delta}) \qquad \mbox{as $\ \eps \rightarrow
0^+$.}
\end{equation*}
\end{thm}

Part (i) of this result was already proved in \cite{RG}.

\bigskip

\section{The $r^{\mathrm{th}}$ moment of the first exit time
for cross-like scatterers}

In this section we estimate the average of the first exit time for
cross-like scatterers, thus proving Theorem \ref{T1.2}. The first
step towards estimating the integral $\int_I l_\eps^r (\omega) \,
d\omega=\int_I l_\eps^{r-2} (\omega) l_\eps^2 (\omega) \, d\omega$
consists in approximating $l_\eps^{r-2}$ by a step function.

We take $I=[\alpha,\beta]\subseteq [ 0,\frac{\pi}{4} ]$, $t_1
=\tan \alpha$, $t_2 =\tan \beta$, $J=[t_1,t_2]\subseteq [0,1]$.
For consecutive Farey fractions
$\frac{a}{q}<\frac{a^\prime}{q^\prime}$ from $J\cap \FF_Q$, where
$Q=\big[ \frac{1}{\eps}\big]$, we denote
\begin{equation*}
\omega_1 =\arctan \frac{a}{q}\, ,\quad \omega_2 =\arctan
\frac{a+\eps}{q} \, ,\quad \omega_2^\prime =\arctan \frac{a^\prime
-\eps}{q^\prime} \, ,\quad \omega_3=\arctan
\frac{a^\prime}{q^\prime} \, .
\end{equation*}

The function $l_\eps^{r-2}$ will be approximated by the constants
$l_\eps^{r-2} (\omega_1) =\| (q,a)\|^{r-2}$ on $[\omega_1,\omega_2
]$ and by $l_\eps^{r-2} (\omega_3)=\| (\q,a^\prime)\|^{r-2}$ on
$[\omega_2,\omega_3 ]$ when $q<\q$, and respectively by
$l_\eps^{r-2} (\omega_1)$ on $[\omega_1,\omega_2^\prime ]$ and by
$l_\eps^{r-2} (\omega_3)$ on $[\omega_2^\prime ,\omega_3 ]$ when
$q>\q$. To be precise, we set
\begin{equation*}
\begin{split}
A_{r,J,\eps} & \hspace{-3pt} =\sideset{^J}{} \sum\limits_{a/q} \|
(q,a)\|^{r-2} \hspace{-4pt} \int\limits_{\omega_1}^{\omega_2}
l_\eps^2 (\omega)\, d\omega + \sideset{^J}{} \sum\limits_{a/q} \|
(\q,a^\prime)\|^{r-2} \hspace{-4pt}
\int\limits_{\omega_2}^{\omega_3} l_\eps^2 (\omega) \, d\omega ,\\
B_{r,J,\eps} & \hspace{-3pt} =\sideset{}{^J} \sum\limits_{a/q}  \|
(q,a)\|^{r-2} \hspace{-4pt}
\int\limits_{\omega_1}^{\omega_2^\prime} l_\eps^2 (\omega) \,
d\omega +\sideset{}{^J} \sum\limits_{a/q} \| (\q,a^\prime)\|^{r-2}
\hspace{-4pt} \int\limits_{\omega_2^\prime}^{\omega_3} l_\eps^2
(\omega) \, d\omega  ,\\ S_{r,J,\eps} & =A_{r,J,\eps}
+B_{r,J,\eps} .
\end{split}
\end{equation*}
Next we estimate the quantities
\begin{equation*}
E_{r,J,\eps}^{(1)} =\Bigg| \sideset{^J}{} \sum\limits_{a/q} \
\int\limits_{\omega_1}^{\omega_3} l_\eps^r (\omega) \, d\omega
-A_{r,J,\eps} \Bigg| \qquad \Big(\leq E^{(1)}_{r,[0,1],\eps}\Big),
\end{equation*}
and respectively
\begin{equation*}
E_{r,J,\eps}^{(2)} =\Bigg| \sideset{}{^J} \sum\limits_{a^\prime
/\q}\ \int\limits_{\omega_1}^{\omega_3} l_\eps^r (\omega) \,
d\omega -B_{r,J,\eps} \Bigg| \qquad \Big(\leq
E^{(2)}_{r,[0,1],\eps}\Big) .
\end{equation*}

An inspection of the case $q<\q$ in the proof of Lemma \ref{L3.1}
leads to
\begin{equation*}
\sup\limits_{\omega \in [\omega_1,\omega_2]} \left| l_\eps^{r-2}
(\omega) -l_\eps^{r-2} (\omega_1) \right| \leq \|
(q,a+\eps)\|^{r-2}-\|(q-\eps,a)\|^{r-2} \ll_r \eps Q^{r-3} ,
\end{equation*}
and to
\begin{equation*}
\begin{split}
\sup\limits_{\omega \in [\omega_2,\omega_3]} & \left| l_\eps^{r-2}
(\omega) -l_\eps^{r-2} (\omega_3) \right| \leq \|
(\q,a^\prime)\|^{r-2} -\left\|
\left( q^\prime,\frac{(a+\eps)q^\prime}{q}\right) \right\|^{r-2} \\
& \ll_r \bigg( a^{\prime } -\frac{(a+\eps) q^{\prime}}{q} \bigg)
Q^{r-3} \ll \frac{Q^{r-3}}{q}\, .
\end{split}
\end{equation*}
Therefore
\begin{equation}\label{4.1}
\begin{split}
\Bigg| \int\limits_{\omega_1}^{\omega_3} l_\eps^r (\omega) \,
d\omega & -l_\eps^{r-2} (\omega_1)
\int\limits_{\omega_1}^{\omega_2} l_\eps^2 (\omega) \, d\omega
-l_\eps^{r-2} (\omega_3) \int\limits_{\omega_2}^{\omega_3}
l_\eps^2 (\omega) \, d\omega \Bigg| \\ & \ll_r q^2 (\omega_2
-\omega_1)\eps Q^{r-3} +q^{\prime 2} (\omega_3 -\omega_2)\,
\frac{Q^{r-3}}{q} \, .
\end{split}
\end{equation}
But $\omega_2 -\omega_1 \leq \frac{a+\eps}{q}-\frac{a}{q}
=\frac{\eps}{q}$ and $\omega_3 -\omega_2 \leq \frac{a^\prime}{\q}
-\frac{a+\eps}{q}=\frac{1-\eps \q}{q\q} <\frac{1}{q\q}$, so the
right-hand side in \eqref{4.1} is
\begin{equation*}
\ll_r q\eps^2 Q^{r-3}+\frac{Q^{r-2}}{q^2} \ll \frac{Q^{r-2}}{q^2}
\, .
\end{equation*}
As a result we infer that
\begin{equation}\label{4.2}
E_{r,[0,1],\eps}^{(1)} =O_r \Bigg(  Q^{r-2} \sum\limits_{q=1}^Q
\frac{\varphi (q)}{q^2} \Bigg) =O_r (Q^{r-2} \ln Q).
\end{equation}

In the case $\q <q$ we get (in both subcases $t_{S^\prime}\leq
t_W$ and $t_{S^\prime}>t_W$)
\begin{equation*}
\begin{split}
\sup\limits_{\omega \in [\omega_1,\omega_2^\prime]} \left|
l_\eps^{r-2}(\omega)-l_\eps^{r-2}(\omega_1)\right| & \leq \left\|
\left( q,\frac{(a^\prime -\eps)q}{q^\prime} \right)
\right\|^{r-2} -\| (q,a-\eps)\|^{r-2} \\
 & \ll_r
\left( \frac{(a^\prime -\eps)q}{q^\prime} -a+\eps \right)
 Q^{r-3} \leq \frac{Q^{r-3}}{q^\prime}
\end{split}
\end{equation*}
and
\begin{equation*}
\sup\limits_{\omega \in [\omega_2^\prime,\omega_3]} \left|
l_\eps^{r-2} (\omega)-l_\eps^{r-2} (\omega_3)\right| \leq \|
(q^\prime,a^\prime)\|^{r-2} -\| (q^\prime,a^\prime -\eps) \|^{r-2}
\ll_r \eps Q^{r-3}.
\end{equation*}

Employing also $\omega_2^\prime -\omega_1 =\frac{a^\prime
-\eps}{\q}-\frac{a}{q} =\frac{1-\eps q}{q\q} \leq \frac{1}{q\q}$
and $\omega_3 -\omega_2^\prime = \frac{a^\prime}{\q} -
\frac{a^\prime -\eps}{\q} =\frac{\eps}{\q}$, we get in the case
$q^\prime <q$ the estimate
\begin{equation*}
\Bigg| \ \int\limits_{\omega_1}^{\omega_3} \hspace{-3pt} l_\eps^r
(\omega) \, d\omega -l_\eps^{r-2} (\omega_1) \hspace{-3pt}
\int\limits_{\omega_1}^{\omega_2^\prime} \hspace{-3pt} l_\eps^2
(\omega) \, d\omega -l_\eps^{r-2} (\omega_3) \hspace{-3pt}
\int\limits_{\omega_2^\prime}^{\omega_3} \hspace{-3pt} l_\eps^2
(\omega) \, d\omega \Bigg|  \ll_r \frac{Q^{r-1}}{qq^{\prime 2}}
+\frac{Q^{r-3}}{q^\prime} \ll \frac{Q^{r-2}}{q^\prime} \, .
\end{equation*}
Hence
\begin{equation}\label{4.3}
E^{(2)}_{r,[0,1],\eps} =O_r \Bigg( \sum\limits_{(q,\q)\in
\Delta_Q} \frac{Q^{r-2}}{q^{\prime 2}} \Bigg) = O_r \Bigg( Q^{r-2}
\sum\limits_{\q =1}^Q \frac{\varphi (\q)}{q^{\prime 2}} \Bigg)
=O_r (Q^{r-2}\ln Q).
\end{equation}
Since the contribution of a single term
$\int_{\omega_1}^{\omega_3} l_\eps^r (\omega) \, d\omega$ is $\ll
\frac{\max \{ q,\q \}^r}{q\q} \leq Q^{r-1}$, we infer from
\eqref{4.2} and \eqref{4.3} that
\begin{equation}\label{4.4}
\int\limits_I l_\eps^r (\omega) \, d\omega =S_{r,J,\eps} +O_r
(\eps^{1-r}).
\end{equation}

Next, we adjust the second integral in the expression of
$A_{r,J,\eps}$, writing
\begin{equation*}
\begin{split}
a^{\prime 2}+q^{\prime 2} & =q^{\prime 2} \left( \left(
\frac{a^\prime}{\q} \right)^{\hspace{-2pt} 2} +1\right) =q^{\prime
2} \left( \left( \frac{a}{q} +\frac{1}{q\q} \right)^{\hspace{-2pt}
2} +1\right) \\ & =\bigg( \frac{\q}{q} \bigg)^{\hspace{-2pt} 2}
\Bigg( \bigg( a+\frac{1}{\q} \bigg)^{\hspace{-2pt} 2} +q^2 \Bigg)
\\ & =\bigg( \frac{\q}{q} \bigg)^{\hspace{-2pt} 2} \Bigg(
a^2+q^2+O\bigg( \frac{a}{\q} \bigg) \Bigg) \\ & =\bigg(
\frac{\q}{q} \bigg)^{\hspace{-2pt} 2} (a^2+q^2) \Bigg(
1+O\bigg( \frac{a}{\q (a^2+q^2)}\bigg) \Bigg) \\
& =\bigg( \frac{\q}{q} \bigg)^{\hspace{-2pt} 2} (a^2 +q^2 )\Bigg(
1+O \bigg( \frac{1}{q\q} \bigg) \Bigg) \\ & =\bigg( \frac{\q}{q}
\bigg)^{\hspace{-2pt} 2} (a^2+q^2)\Bigg( 1+O \bigg( \frac{1}{Q}
\bigg) \Bigg) .
\end{split}
\end{equation*}

This gives in turn
\begin{equation}\label{4.5}
\| (\q,a^\prime)\|^{r-2} =\bigg( \frac{\q}{q} \bigg)^{r-2}
\|(q,a)\|^{r-2} \bigg( 1+O_r \Bigl( \frac{1}{Q} \Bigr) \bigg) .
\end{equation}

In a similar way
\begin{equation}\label{4.6}
\|(q,a)\|^{r-2} =\bigg( \frac{q}{\q} \bigg)^{\hspace{-2pt} r-2}
\|(\q,a^\prime)\|^{r-2} \bigg( 1+O_r \Bigl( \frac{1}{Q} \Bigr)
\bigg) .
\end{equation}

For further use, it is also worth to note
\begin{equation}\label{4.7}
\frac{1+\frac{a^\prime}{\q}}{1+(\frac{a^\prime}{\q})^2}
=\frac{\big( 1+\frac{a}{q} \big) \big( 1+O ( \frac{1}{Q}
)\big)}{\big( 1+(\frac{a}{q})^2 \big) \big( 1+O( \frac{1}{Q} )
\big)} =\frac{1+\frac{a}{q}}{1+(\frac{a}{q})^2} \, \Bigg( 1+O
\bigg( \frac{1}{Q} \bigg) \Bigg) .
\end{equation}

Making use of (see the proof of Lemma \ref{L3.1})
\begin{equation*}
\begin{split}
\int\limits_{\omega_1}^{\omega_2} l_\eps^2 (\omega) \, d\omega &
=2\ar (\triangle OAN) +O (\eps^2) =\eps q+O (\eps^2),
\\ \int\limits_{\omega_2}^{\omega_3} l_\eps^2 (\omega) \, d\omega
& =2\ar (\triangle ON_0^\prime A^\prime)=\frac{\q (1-\eps \q)}{q}
\, ,
\end{split}
\end{equation*}
and \eqref{4.5}, we see that $A_{r,J,\eps}$ can be expressed as
\begin{equation}\label{4.8}
\begin{split}
\sideset{^J}{} \sum\limits_{a/q} \|(q,a)\|^{r-2} \bigg(
\frac{\q}{q} \bigg)^{\hspace{-2pt} r-2} & \frac{\q(1-\eps \q)}{q}
\bigg( 1+O_r \Bigl( \frac{1}{Q} \Bigr) \bigg) \\ & +\eps
\sideset{^J}{} \sum\limits_{a/q} \|(q,a)\|^{r-2} q +O_r (1) .
\end{split}
\end{equation}

If $\X$ denotes the inverse of the integer $x\hspace{-2pt}
\pmod{q}$ in $[1,q]$, then $a^\prime q-a\q =1$ gives
$a=q-\bar{q^\prime}$. Since $\frac{1-\eps \q}{q}<\eps$, the error
in the first sum in \eqref{4.8} is $\ll_r Q^2 Q^{r-1} \eps Q^{-1}
=Q^{r-1}$, and so $A_{r,J,\eps}$ is equal up to an error term of
order $O_r (Q^{r-1})$ to
\begin{equation}\label{4.9}
\begin{split}
& \sideset{^J}{} \sum\limits_{a/q} \|(q,a)\|^{r-2} \bigg( \eps
q+\frac{q^{\prime r-1}}{q^{r-1}} -\frac{\eps q^{\prime
r}}{q^{r-1}} \bigg) \\ & =
\sum\limits_{q=1}^Q \ \sum_{\substack{\max \{ q,Q-q\}<x\leq Q \\
\X \in J_q^{(1)}}} \hspace{-5pt} \big( (q-\X)^2 + q^2
\big)^{\frac{r-2}{2}} \bigg( \eps q+\frac{x^{r-1}}{q^{r-1}}
-\frac{\eps x^r}{q^{r-1}} \bigg) ,
\end{split}
\end{equation}
with $J_q^{(1)}$ as defined in Section 2.

When $q>\q$, we employ (see the proof of Lemma \ref{L3.1})
\begin{equation*}
\begin{split}
\int\limits_{\omega_1}^{\omega_2^\prime} l_\eps^2 (\omega) \,
d\omega & =2\ar (\triangle OAS_0^\prime) +O
(\eps^2)=\frac{(1-\eps q)q}{\q} +O (\eps^2) , \\
\int\limits_{\omega_2^\prime}^{\omega_3} l_\eps^2 (\omega) \,
d\omega & =2\ar (\triangle OA^\prime S^\prime) =\eps \q ,
\end{split}
\end{equation*}
together with $Q^{r-2} \eps^2 \# \FF_Q \leq Q^{r-2}$, relation
\eqref{4.6}, and the fact that $a^\prime q-a\q =1$ implies
$a^\prime =\bar{q}\hspace{-2pt} \pmod{\q}$, to infer that
$B_{r,J,\eps}$ is expressible as
\begin{equation*}
\begin{split}
& \sideset{}{^J} \sum\limits_{a^\prime /\q}  \|(q,a)\|^{r-2}
\frac{(1-\eps q)q}{\q} + \sideset{}{^J} \sum\limits_{a^\prime /\q}
\hspace{-3pt} \bigg( \frac{\q}{q} \bigg)^{\hspace{-3pt} r-2}
\hspace{-10pt} \|(q,a)\|^{r-2} \eps \q +O_r (Q^{r-1})
\\ & =\sideset{}{^J} \sum\limits_{a^\prime /\q} \|
(\q,a^\prime)\|^{r-2} \bigg( \eps \q +\frac{q^{r-1}}{q^{\prime
r-1}} -\frac{\eps q^r}{q^{\prime r-1}} \bigg) +O_r (Q^{r-1})
\\ & = \sum\limits_{\q =1}^Q \, \sum_{\substack{\max \{ \q
,Q-\q \} <x\leq Q \\ \X \in J_{\q}^{(2)}}} \hspace{-10pt}
\hspace{-15pt}\big( \X^2 +q^{\prime 2} \big)^{\frac{r-2}{2}}
\bigg( \eps \q +\frac{x^{r-1}}{q^{\prime r-1}} -\frac{\eps
x^r}{q^{\prime r-1}} \bigg) +O_r (Q^{r-1}) ,
\end{split}
\end{equation*}
where $\X$ denotes the inverse of an integer $x\hspace{-2pt}
\pmod{\q}$ in $[1,\q]$. Changing notation, $B_{r,J,\eps}$ can be
rewritten as
\begin{equation}\label{4.10}
\sum\limits_{q=1}^Q \ \sum_{\substack{\max \{ q,Q-q\}<x\leq Q
\\ \X \in J_q^{(2)}}} \hspace{-9pt} \big( \X^2+q^2
\big)^{\frac{r-2}{2}} \bigg( \eps q+\frac{x^{r-1}}{q^{r-1}}
-\frac{\eps x^r}{q^{r-1}} \bigg) +O_r (Q^{r-1}).
\end{equation}

By \eqref{4.4}, \eqref{4.9} and \eqref{4.10}, we infer that
\begin{equation}\label{4.11}
\int\limits_I l_\eps^r (\omega) \, d\omega = T_{r,J,\eps}
+O_r(Q^{r-1}),
\end{equation}
where $T_{r,J} (\eps)=S_1(Q)+S_2(Q)$, with
\begin{equation}\label{4.12}
S_k(Q) =\sum\limits_{q=1}^Q \, \sum_{\substack{\max \{
Q-q,q\}<x\leq Q \\ \X \in J_q^{(k)}}} f_k (x,\X,q),\qquad k=1,2,
\end{equation}
and
\begin{equation*}
\begin{split}
& f_2 (x,y,z) = \big( y^2+z^2 \big)^{\frac{r-2}{2}} \bigg( \eps
z+\frac{x^{r-1}}{z^{r-1}} -\frac{\eps x^r}{z^{r-1}}
\bigg) , \\
&  f_1 (x,y,z)=f_2 (x,z-y,z).
\end{split}
\end{equation*}

For each $q\in [1,Q]$, the functions $f_k (\cdot ,\cdot, q)$,
defined on $[1,Q]\times [1,q]$, manifestly satisfy the estimates
\begin{equation}\label{4.13}
\| f_k (\cdot ,\cdot ,q)\|_{\infty} \ll_r \frac{Q^{r-1}}{q}
\end{equation}
and
\begin{equation}\label{4.14}
\| Df_k (\cdot ,\cdot ,q)\|_{\infty} \ll_r \frac{Q^{r-2}}{q} \leq
\frac{Q^{r-1}}{q^2} \, .
\end{equation}
Thus we may consider in Lemma \ref{L2.2} for each $q\in [1,Q]$ the
function $f_k (\cdot ,\cdot ,q)$, the intervals $\II =( \max
\{Q-q,q\} ,Q]$ and $\JJ =J_q^{(k)}$ with $\vert \II \vert \leq q$
and $\vert \JJ \vert =q\vert I\vert$, and take $T=[Q^\alpha ]$, to
infer that the inner sum in \eqref{4.12} can be expressed as
\begin{equation*}
\begin{split}
\frac{\varphi (q)}{q^2} & \int\limits_{\max \{Q-q,q\}}^Q
\hspace{-8pt}
dx \int\limits_{J_q^{(k)}} dy\ f_k (x,y,q) \\
&  +O_{\delta,r} ( q^{-\frac{1}{2}+\delta} Q^{r-1+2\alpha}
+q^{-\frac{1}{2}+\delta} Q^{r-1+\alpha} +\vert I\vert
Q^{r-1-\alpha} ) .
\end{split}
\end{equation*}

Summing up over $q\in [1,Q]$, we arrive at
\begin{equation}\label{4.15}
S_k(Q) =\sum\limits_{q=1}^Q \frac{\varphi (q)}{q} \, g(q)+
O_{\delta,r} ( Q^{r-\frac{1}{2}+2\alpha+\delta} +\vert I\vert
Q^{r-\alpha} ) ,
\end{equation}
where
\begin{equation*}
\begin{split}
g (z) & =\frac{1}{z} \int\limits_{\max \{ Q-z,z\}}^Q
\hspace{-10pt} dx \hspace{3pt} \int\limits_{(1-t_2)z}^{(1-t_1)z}
\hspace{-4pt} dy\  f_1 (x,y,z)  \\  & =\frac{1}{z}
\int\limits_{\max \{Q-z,z\}}^Q \hspace{-10pt} dx \hspace{3pt}
\int\limits_{t_1 z}^{t_2 z} dy\  f_2 (x,y,z),\qquad z\in [1,Q].
\end{split}
\end{equation*}
The formulas
\begin{equation*}
\begin{split}
& \frac{d}{dz} \Bigg( \frac{1}{z} \int\limits_z^Q dx
\int\limits_{az}^{bz} dy\, h(x,y,z)\Bigg) = -\frac{1}{z^2}
\int\limits_z^Q dx \int\limits_{az}^{bz} dy\, h(x,y,z) \\ &
\hspace{3.5cm} +\frac{1}{z} \int\limits_z^Q dx
\int\limits_{az}^{bz} dy\, \frac{\partial h}{\partial z} \,
(x,y,z)-\frac{1}{z}
\int\limits_{az}^{bz} h(z,y,z)\, dy \\
& \hspace{3.5cm} +\frac{b}{z} \int\limits_z^Q
h(x,bz,z)\, dx-\frac{a}{z} \int\limits_z^Q h(x,az,z)\, dx,\\
& \frac{d}{dz} \Bigg( \frac{1}{z} \int\limits_{Q-z}^Q dx
\int\limits_{az}^{bz} dy\, h(x,y,z) \Bigg) = -\frac{1}{z^2}
\int\limits_{Q-z}^Q dx \int\limits_{az}^{bz} dy\, h(x,y,z) \\
& \hspace{3.5cm} +\frac{1}{z} \int\limits_{Q-z}^Q dx
\int\limits_{az}^{bz} dy\, \frac{\partial h}{\partial z} \,
(x,y,z)  +\frac{1}{z} \int\limits_{az}^{bz} h(Q-z,y,z)\, dy \\ &
\hspace{3.5cm} +\frac{b}{z} \int\limits_{Q-z}^Q h(x,bz,z)\,
dx-\frac{a}{z} \int\limits_{Q-z}^Q h(x,az,z)\, dx,
\end{split}
\end{equation*}
and the estimates \eqref{4.13} and \eqref{4.14} show that $\vert
g^\prime (z)\vert \ll_r \frac{Q^{r-1}}{z}$. As a result we get
$\int_1^Q \vert g^\prime (z)\vert \, dz \ll_r Q^{r-1} \ln Q$. It
is also clear that $\| g \|_\infty \ll_r Q^{r-1}$, so we are in
the position of being able to apply Lemma 2.3 in \cite{BCZ1} to
$g$, collecting
\begin{equation*}
\begin{split}
\sum\limits_{q=1}^Q \frac{\varphi (q)}{q}\, g (q) &
=\frac{1}{\zeta (2)} \int\limits_1^Q g (z)\, dz+O\Bigg( \bigg( \|
g \|_\infty +\int\limits_1^Q \vert g^\prime (z)\vert \, dz \bigg)
\ln Q \Bigg) \\ & =\frac{1}{\zeta (2)} \int\limits_1^Q g (z)\, dz
+O_r (Q^{r-1} \ln^2 Q) .
\end{split}
\end{equation*}

Comparing the previous relation with \eqref{4.15}, we infer that
both $S_1(Q)$ and $S_2(Q)$ can now be expressed as
\begin{equation}\label{4.16}
\frac{1}{\zeta(2)} \int\limits_1^Q g(z)\, dz +O_{\delta,r} (
Q^{r-\frac{1}{2}+2\alpha+\delta}+\vert I\vert Q^{r-\alpha} ) .
\end{equation}

Taking into account \eqref{4.11} and \eqref{4.12}, we gather
\begin{equation}\label{4.17}
\int\limits_I l_\eps^r (\omega) \, d\omega =\frac{2}{\zeta(2)}
\int\limits_1^Q g(z)\, dz +O_{\delta,r}  (
Q^{r-\frac{1}{2}+2\alpha+\delta}+\vert I\vert Q^{r-\alpha} ) .
\end{equation}

Integrating with respect to $y$ in the formula that gives $g$ and
changing then $z$ into $Qz$, we may express the main term in
\eqref{4.17} as
\begin{equation}\label{4.18}
\begin{split}
& \frac{K_{r,I}}{\zeta(2)} \int\limits_1^Q dz \int\limits_{\max
\{Q-z,z\}}^Q \hspace{-8pt} dx \ \bigg( \eps
+\frac{x^{r-1}}{z^r} -\frac{\eps x^r}{z^r} \bigg) z^{r-1} \\
& =\frac{K_{r,I} Q^{r+1}}{\zeta(2)} \int\limits_{1/Q}^1 dz
\int\limits_{\max \{1-z,z\}}^1 \hspace{-8pt} dx\ \bigg( \eps
z^{r-1} +\frac{x^{r-1}}{Qz} -\frac{\eps x^r}{z} \bigg),
\end{split}
\end{equation}
where
\begin{equation*}
K_{r,I} =2\int\limits_{t_1}^{t_2} (1+t^2)^{\frac{r-2}{2}} dt =2
\int\limits_I \frac{dx}{\cos^r x} \, .
\end{equation*}
Up to an error term of order $O_r(Q^{-2})$, the double integral in
the right-hand side of \eqref{4.18} is given by
\begin{equation*}
\varepsilon \hspace{-1.5pt} \int\limits_0^{1/2} \hspace{-2pt}
z(z^{r-1}+(1-z)^{r-1})\, dz +\frac{1}{rQ} \int\limits_0^{1/2}
\frac{1-(1-z)^r}{z(1-z)} \, dz-\frac{\varepsilon}{r+1}
\int\limits_0^{1/2} \frac{1-(1-z)^{r+1}}{z(1-z)} \, dz .
\end{equation*}
Since $Q^{r+1}=\frac{1}{\eps^{r+1}} +O_r ( \frac{1}{\eps^r})$, we
now infer from \eqref{4.17} and \eqref{4.18} the equality
\begin{equation}\label{4.19}
\int\limits_I l_\eps^r (\omega) \, d\omega = \frac{c_r}{\eps^r}
\int\limits_I \frac{dx}{\cos^r x} +O_{\delta,r}
(\eps^{-r+\frac{1}{2}-2\alpha-\delta}+\vert I\vert
\eps^{-r+\alpha}),
\end{equation}
with $c_r$ as in Theorem \ref{T1.1}. The proof of Theorem
\ref{T1.2} is now complete.

\bigskip

\section{The moments of the number of reflections}

We take as before $Q=[\frac{1}{\eps}]$ and keep up with the
notation from the beginning of Section 4. An inspection of the
proof of Lemma \ref{L3.1} shows that if $\frac{a}{q}
<\frac{a^\prime}{q^\prime}$ are consecutive in $\FF_Q$, then
\begin{equation*}
R_\eps (\omega)=\begin{cases} & \begin{cases} q+a & \mbox{\rm if
$\tan \omega \in [\frac{a}{q},
\frac{a}{q-\eps} )$} \\
q+a+1 & \mbox{\rm if $\tan \omega \in [\frac{a}{q-\eps} ,
\frac{a+\eps}{q} )$} \\
q^\prime +a^\prime & \mbox{\rm if $\tan \omega \in
[\frac{a+\eps}{q}, \frac{a^\prime}{q^\prime} ]$}
\end{cases} \qquad \mbox{\rm if $q<q^\prime$;} \\
& \begin{cases} q+a & \mbox{\rm if $\tan \omega \in [\frac{a}{q},
\frac{a^\prime -\eps}{q^\prime})$} \\
q^\prime +a^\prime & \mbox{\rm if $\tan \omega \in [\frac{a^\prime
-\eps}{q^\prime}, \frac{a^\prime}{q^\prime}]$,}
\end{cases} \qquad \qquad \mbox{\rm if $q>q^\prime$ and
$t_{S^\prime} \leq t_W$;} \\
& \begin{cases}
q+a & \mbox{\rm if $\tan \omega \in [\frac{a}{q},\frac{a}{q-\eps})$}\\
 q+a+1 & \mbox{\rm if $\tan \omega \in [\frac{a}{q-\eps},
\frac{a^\prime -\eps}{q^\prime})$} \\
q^\prime +a^\prime & \mbox{\rm if $\tan \omega \in [\frac{a^\prime
-\eps}{q^\prime} ,\frac{a^\prime}{q^\prime}],$}
\end{cases} \qquad \mbox{\rm if $q>q^\prime$ and
$t_{S^\prime} >t_W$.}
\end{cases}
\end{equation*}

A first immediate remark is that we may replace $q+a+1$ by $q+a$
in the above formulas, since the contribution of the corresponding
arcs is small, as we see from $\vert \arctan x-\arctan y\vert \leq
\vert x-y\vert$, and from
\begin{equation*}
\sum\limits_{a/q\in \FF_Q} \left( \frac{a+\eps}{q} -
\frac{a}{q-\eps} \right) =\sum\limits_{a/q \in \FF_Q} \frac{\eps
(q-a-\eps)}{q(q-\eps)} \ll \eps \sum\limits_{a/q\in \FF_Q}
\frac{1}{q} \leq 1
\end{equation*}
and
\begin{equation*}
\sum\limits_{a/q\in \FF_Q} \left( \frac{a^\prime -\eps}{q^\prime}
-\frac{a}{q-\eps} \right) =\sum\limits_{a/q\in \FF_Q} \frac{1-\eps
(q+a^\prime)+\eps^2}{q^\prime (q-\eps)} \ll \sum\limits_{a/q \in
\FF_Q} \frac{1}{qq^\prime} =1.
\end{equation*}
As a result we may write
\begin{equation*}
\int\limits_I R_\eps^r (\omega)\, d\omega =T^{(1)}_{r,J,\eps}
+T^{(2)}_{r,J,\eps} +O_r (Q^{r-1}) ,
\end{equation*}
where we set
\begin{equation*}
\begin{split}
T^{(1)}_{r,J,\eps} =\sideset{^J}{} \sum\limits_{a/q} & (q+a)^r
\left( \arctan \frac{a+\eps}{q} -\arctan \frac{a}{q}
\right)  \\
& \qquad + \sideset{^J}{} \sum\limits_{a/q} (q^\prime +a^\prime)^r
\left( \arctan \frac{a^\prime}{q^\prime} -\arctan \frac{a+\eps}{q}
\right)
\end{split}
\end{equation*}
and
\begin{equation*}
\begin{split}
T^{(2)}_{r,J,\eps} =\sideset{}{^J} \sum\limits_{a/q} & (q+a)^r
\left( \arctan \frac{a^\prime -\eps}{q^\prime} -
\arctan \frac{a}{q}\right)  \\
& \qquad + \sideset{}{^J} \sum\limits_{a/q} (q^\prime +a^\prime)^r
\left( \arctan \frac{a^\prime}{q^\prime} -\arctan \frac{a^\prime
-\eps}{q^\prime} \right).
\end{split}
\end{equation*}

Employing
\begin{equation*}
\arctan (x+h)-\arctan x=\frac{h}{1+x^2} +O(h^2)=
\frac{h}{1+(x+h)^2} +O(h^2)
\end{equation*}
we now arrive at
\begin{equation}\label{5.1}
\int\limits_I R_\eps^r (\omega) \, d\omega =
S_{r,J,\eps}+T_{r,J,\eps} +O_r (Q^{r-1} ),
\end{equation}
where
\begin{equation}\label{5.2}
S_{r,J,\eps} =\sideset{^J}{} \sum\limits_{a/q}
 \left( \frac{\frac{\eps}{q}}{1+
(\frac{a}{q})^2} \, (q+a)^r+ \frac{\frac{1-\eps
q^\prime}{qq^\prime}}{1+ (\frac{a^\prime}{q^\prime})^2} \,
(q^\prime+a^\prime)^r \right)
\end{equation}
and
\begin{equation}\label{5.3}
T_{r,J,\eps} =\sideset{}{^J} \sum\limits_{a/q}  \left(
\frac{\frac{1-\eps q}{qq^\prime}}{1+(\frac{a}{q})^2} \, (q+a)^r
+\frac{\frac{\eps}{q^\prime}}{1+(\frac{a^\prime}{q^\prime})^2} \,
(q^\prime+a^\prime)^r \right) .
\end{equation}
To further simplify the expressions in \eqref{5.2} and \eqref{5.3}
we employ
\begin{equation*}
\begin{split}
& 1+\frac{a^\prime}{q^\prime} =1+\frac{a}{q}+\frac{1}{qq^\prime}
=1+\frac{a}{q}+O\left( \frac{1}{Q} \right)= \left( 1+\frac{a}{q}
\right) \left( 1+O\Big( \frac{1}{Q}
\Big) \right) ,\\
& \left( 1+\frac{a^\prime}{q^\prime} \right)^r =\left( 1+
\frac{a}{q} \right)^r \left( 1+O_r \Big( \frac{1}{Q} \Big) \right)
,
\end{split}
\end{equation*}
and
\begin{equation*}
\frac{1}{Q} \sum\limits_{a/q\in \FF_Q} \frac{1}{qq^\prime}\,
q^{\prime r} \leq Q^{r-2} \sum\limits_{a/q\in \FF_Q} \frac{1}{q}
\leq Q^{r-1},
\end{equation*}
to infer that
\begin{equation}\label{5.4}
S_{r,J,\eps}=\sideset{^J}{} \sum\limits_{a/q} \frac{
(1+\frac{a}{q})^r}{1+(\frac{a}{q})^2} \left( \eps
q^{r-1}+\frac{1-\eps q^\prime}{q} \, q^{\prime r-1} \right) +O_r
(Q^{r-1}),
\end{equation}
and also that
\begin{equation*}
T_{r,J,\eps} =\sideset{}{^J} \sum\limits_{a/q} \frac{
(1+\frac{a}{q})^r}{1+(\frac{a}{q})^2} \left( \eps
q^{r-1}+\frac{1-\eps q^\prime}{q} \, q^{\prime r-1} \right) +O_r
(Q^{r-1}).
\end{equation*}
The main term in \eqref{5.4} can now be conveniently expressed as
\begin{equation}\label{5.5}
A_{r,J,\eps}=\sum\limits_{q=1}^Q \sum_{\substack{\max
\{Q-q,q\}<x\leq Q \\ \bar{x} \in J_q^{(1)}}} f(x,\bar{x},q),
\end{equation}
where $\bar{x}$ is the multiplicative inverse of $x$ in $\{
1,\dots,q-1 \}$, and this time we set
\begin{equation*}
\begin{split}
f(x,y) =f(x,y,q) &
=\frac{(1+\frac{q-y}{q})^r}{q^2+(\frac{q-y}{q})^2}
\left( \eps q^{r-1} +\frac{1-\eps x}{q} \, x^{r-1} \right) \\
& =\frac{q(2q-y)^r}{q^2 +(q-y)^2} \left( \eps +\frac{1-\eps x}{q}
\Big( \frac{x}{q}\Big)^{r-1} \right) .
\end{split}
\end{equation*}

Since $0<1-\eps x<1$, it is easy to check that
\begin{equation*}
\begin{split}
\| f\|_{\infty,\II\times J_q^{(1)}} & \ll
\frac{(2q-y)^r}{q^2+(q-y)^2} \ll_r q^{r-2} ,\\
\left\| \frac{\partial f}{\partial x} \right\|_{\infty,\II \times
J_q^{(1)}} & \ll \frac{q(2q-y)^r}{q^2+(q-y)^2} \, \frac{1}{q^2}
\ll_r q^{r-3} ,\\
\left\| \frac{\partial f}{\partial y} \right\|_{\infty , \II\times
J_q^{(2)}} & \ll \left\| \frac{\partial}{\partial y} \left(
\frac{(2q-y)^r}{q^2+(q-y)^2} \right) \right\|_\infty \ll_r q^{r-3}
.
\end{split}
\end{equation*}

Taking these estimates into account and applying Lemma \ref{L2.2}
for each $q\in [1,Q]$ to $f=f(\cdot,\cdot,q)$, $\II=(\max \{
Q-q,q\},Q]$, $\JJ=J_q^{(1)}$ and $T=[Q^\alpha ]$, we infer that
the inner sum in \eqref{5.5} can be expressed as
\begin{equation*}
\begin{split}
\frac{\varphi (q)}{q^2} & \int\limits_{\max \{ Q-q,q\}}^Q \left(
\eps+ \frac{1-\eps x}{q} \Big( \frac{x}{q} \Big)^{r-1} \right) dx
\  \int\limits_{J_q^{(1)}}
\frac{q(2q-y)^r}{q^2+(q-y)^2} \, dy \\
& \qquad +O_{r,\delta} ( Q^{2\alpha} q^{\frac{1}{2}+\delta}
q^{r-2} +Q^\alpha q^{\frac{3}{2}+\delta} q^{r-3} +
\vert I\vert Q^{-\alpha} q^2 q^{r-3} ) \\
& =\frac{C_{r,I} \varphi (q)}{q} \ q^{r-1} g_r (q) +O_{r,\delta} (
Q^{r-\frac{3}{2}+2\alpha +\delta} +\vert I\vert Q^{r-1-\alpha} ),
\end{split}
\end{equation*}
where
\begin{equation*}
\begin{split}
C_{r,I} & =\frac{1}{q^{r-1}} \int\limits_{(1-t_2)q}^{(1-t_1)q}
\frac{(2q-y)^r}{q^2+(q-y)^2} \, dy=\frac{1}{q^{r-1}}
\int\limits_{t_1 q}^{t_2 q}
\frac{(q+y)^r}{q^2+y^2} \, dy \\
& =\int\limits_{t_1}^{t_2} \frac{(1+t)^r}{1+t^2} \, dt
=\int\limits_I (1+\tan x)^r \, dx
\end{split}
\end{equation*}
and
\begin{equation*}
g_r (q)=\int\limits_{\max \{Q-q,q\}}^Q \left( \eps + \frac{1-\eps
x}{q} \Big( \frac{x}{q} \Big)^{r-1} \right) dx.
\end{equation*}

We now arrive at
\begin{equation}\label{5.6}
\sum_{\substack{\max \{Q-q,q\}<x\leq Q \\ \bar{x} \in J_q^{(1)}}}
f(x,\bar{x},q) =\frac{C_{r,I}}{q} \, h_r (q)+ O_{r,\delta} (
Q^{r-\frac{3}{2}+2\alpha+\delta} +\vert I\vert Q^{r-\alpha-1} ) ,
\end{equation}
where
\begin{equation*}
h_r(q)=\int\limits_{\max \{Q-q,q\}}^Q \left( \eps
q^{r-1}+\frac{1-\eps x}{q} \, x^{r-1} \right) dx.
\end{equation*}

Since $\| h\|_\infty \ll Q^{r-1}$ and $\int_1^Q \vert h^\prime
(q)\vert \, dq \ll Q^{r-1}$, Lemma 2.3 in \cite{BCZ1} together
with \eqref{5.5} and \eqref{5.6} show that
\begin{equation*}
\begin{split}
S_{r,J,\eps} & =C_{r,I} \sum\limits_{q=1}^Q \frac{\varphi (q)}{q}
\, h_r (q)+O_{r,\delta} (
Q^{r-\frac{1}{2}+2\alpha+\delta}+\vert I\vert Q^{r-\alpha} ) \\
& =\frac{C_{r,I}}{\zeta (2)} \int\limits_1^Q h_r(q)\,
dq+O_{r,\delta} ( Q^{r-\frac{1}{2}+2\alpha+\delta} +\vert I\vert
Q^{r-\alpha} ) .
\end{split}
\end{equation*}

Making use of $\eps Q=1+O(\eps )$ we arrive by a straightforward
computation to
\begin{equation*}
\begin{split}
\int\limits_1^Q h_r(q)\, dq & =Q^r \int\limits_0^1 \left( x^{r-1}
\big( 1-\max \{1-x,x\} \big) + \frac{1-\max \{1-x,x\}^r}{rx}
\right. \\ & \left. \qquad \qquad \qquad \qquad -
\frac{1-\max\{1-x,x\}^{r+1}}{(r+1)x} \right) dx+ O_r (Q^{r-1}).
\end{split}
\end{equation*}
The integral above is seen to coincide with
\begin{equation*}
\int\limits_0^{1/2} \left( x\big( x^{r-1}+(1-x)^{r-1}\big)
+\frac{1-(1-x)^r}{rx(1-x)}- \frac{1-(1-x)^{r+1}}{(r+1)x(1-x)}
\right) dx =\frac{\pi^2 c_r}{12} \ ,
\end{equation*}
hence
\begin{equation*}
\int\limits_1^Q h_r(q)\, dq=\frac{\pi^2 c_r Q^r}{12}+ O_r(Q^{r-1})
=\frac{\pi^2 c_r \eps^{-r}}{12} +O_r (\eps^{-r+1})
\end{equation*}
and as a result
\begin{equation*}
\begin{split}
S_{r,J,\eps} & =\frac{6}{\pi^2}\cdot \frac{\pi^2 c_r
\eps^{-r}}{12} \int\limits_I (1+\tan x)^r \, dx + O_{r,\delta}
(\eps^{-r+\frac{1}{2}-2\alpha-\delta}
+\vert I\vert \eps^{-r+\alpha}) \\
& =\frac{c_r \eps^{-r}}{2} \int\limits_I (1+\tan x)^r \, dx +
O_{r,\delta} (\eps^{-r+\frac{1}{2}-2\alpha-\delta} +\vert I\vert
\eps^{-r+\alpha}).
\end{split}
\end{equation*}

By reversing the roles of $q$ and $q^\prime$ it is seen in a
similar way that
\begin{equation*}
T_{r,J,\eps} =\frac{c_r \eps^{-r}}{2} \int\limits_I (1+\tan x)^r\,
dx +O_r (\eps^{-r+\frac{1}{2}-2\alpha-\delta} +\vert I\vert
\eps^{-r+\alpha} ).
\end{equation*}
This concludes the estimates of $S_{r,J,\eps}$ and $T_{r,J,\eps}$.
Theorem \ref{T1.5} now follows from \eqref{5.1}.

\bigskip

\section{The case of circular scatterers}

Note first that the statements of Theorems \ref{T1.2} and
\ref{T1.5} hold true if we replace the scatterers $C_\eps
+\Z^{2\ast}$ by $V_\eps +\Z^{2\ast}$.

In this section we consider the circular scatterers
$D_\eps+\Z^{2\ast}$, where
\begin{equation*}
D_\eps =\{(x,y)\in \R^2 \, ;\, x^2+y^2=\eps^2 \} .
\end{equation*}

For each integer lattice point $(q,a)$, let $(q,a\pm \eps_\pm )$
denote the intersections of the line $x=q$ with the tangents from
$O$ to the circle
\begin{equation*}
D_{\eps,q,a}=(q,a)+D_\eps=\{ (x,y)\in \R^2 \,
;\,(x-q)^2+(y-a)^2=\eps^2 \},
\end{equation*}
where $\eps_\pm =\eps_\pm (q,a)$ are computed from the equality
\begin{equation*}
\eps =\frac{\big| a-\frac{a\pm \eps_\pm}{q} q\big|}{\sqrt{1+\big(
\frac{a\pm \eps_\pm}{q} \big)^2}}=\frac{\eps_\pm
q}{\sqrt{q^2+(a\pm \eps_\pm)^2}} \, ,
\end{equation*}
which gives in turn
\begin{equation*}
\eps_\pm^2 q^2=\eps^2 q^2+\eps^2 (a\pm \eps_\pm)^2 ,
\end{equation*}
or
\begin{equation*}
(q^2-\eps^2)\eps_\pm^2 \mp 2a\eps^2 \eps_\pm -\eps^2 (q^2+a^2)=0.
\end{equation*}

\begin{figure}[ht]
\includegraphics*[scale=0.7, bb=0 0 250 220]{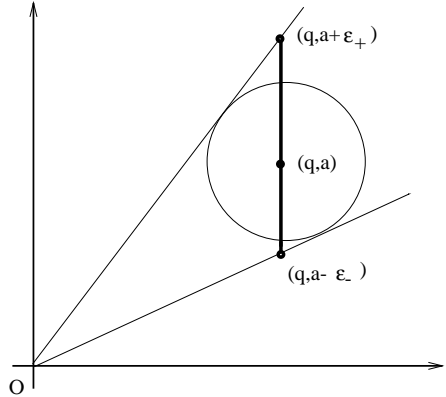}
\caption{A circular scatterer} \label{Figure5}
\end{figure}

The latter provides
\begin{equation}\label{6.1}
\eps_\pm =\pm \frac{a\eps^2}{q^2-\eps^2} +\frac{\eps}{q^2-\eps^2}
\sqrt{q^4+a^2 q^2 -a^2 \eps^2} .
\end{equation}
Employing also
\begin{equation*}
\begin{split}
& \sqrt{q^4+a^2 q^2} -\sqrt{q^4+a^2 q^2 -a^2 \eps^2} \ll \eps^2
,\\
& \frac{\eps}{q^2-\eps^2} =\frac{\eps}{q^2} \bigg( 1+O\Big(
\frac{\eps^2}{q^2} \Big) \bigg)
\end{split}
\end{equation*}
and
\begin{equation*}
\frac{a\eps^2}{q^2-\eps^2} =O\bigg( \frac{\eps^2}{q} \bigg) ,
\end{equation*}
we arrive at
\begin{equation}\label{6.2}
\eps_\pm (q,a)=\eps \sqrt{1+\frac{a^2}{q^2}} +O\bigg(
\frac{\eps^2}{q} \bigg)=\frac{\eps}{\cos \arctan \frac{a}{q}}
+O\left( \frac{\eps^2}{q} \right) .
\end{equation}

\begin{proof}[Proof of Theorem \ref{T1.1}]
We wish to compare $\int_I \tilde{\tau}_\eps^r (\omega)\, d\omega$
with $\int_I \tilde{\tilde{\tau}}_\eps^r (\omega)\, d\omega$ where
$\tilde{\tilde{\tau}}_\eps (\omega)$, the smallest $\tau >0$ for
which
\begin{equation*}
(\tau \cos \omega,\tau \sin \omega )\in \bigcup\limits_{(q,a)\in
\Z^{2 \ast}} \hspace{-5pt} \{ q\} \times [a-\eps_-
(q,a),a+\eps_+(q,a)] ,
\end{equation*}
denotes the first exit time in the case where the scatterers are
the vertical segments $\{ q\} \times [a-\eps_-(q,a),
a+\eps_+(q,a)]$. From Figure \ref{Figure5} it is apparent that
\begin{equation*}
\sup\limits_\omega \vert \tilde{\tau}_\eps (\omega)-
\tilde{\tilde{\tau}}_\eps (\omega)\vert \leq 2\eps ,
\end{equation*}
and so, since $\sup_\omega \tilde{\tilde{\tau}}_\eps (\omega) \leq
\sup_\omega l_\eps (\omega) \leq \frac{\sqrt{2}}{\eps}$, we get
\begin{equation*}
\sup\limits_\omega \vert \tilde{\tau}_\eps^r (\omega) -
\tilde{\tilde{\tau}}_\eps^r (\omega) \vert \ll_r \eps \bigg(
\frac{\sqrt{2}}{\eps} \bigg)^{r-1} \ll_r \eps^{2-r} ,
\end{equation*}
which gives
\begin{equation}\label{42}
\int\limits_I \tilde{\tau}_\eps^r (\omega)\, d\omega =
\int\limits_I \tilde{\tilde{\tau}}_\eps^r (\omega)\, d\omega +O_r
(\eps^{2-r}).
\end{equation}

To estimate $\int_I \tilde{\tilde{\tau}}_\eps (\omega)\, d\omega$,
we divide the interval $I$ into $N=[\eps^{-\theta}]$ intervals of
equal size $I_j =[\omega_j,\omega_{j+1}]$ with $\vert I_j \vert
=\frac{\vert I\vert}{N}\asymp \eps^\theta$ for some $0<\theta
<\frac{1}{2}$. Then one has for all $j$ that
\begin{equation*}
\left| \frac{\cos \omega_j}{\eps-\eps^{3/2}} - \frac{\cos
\omega_{j+1}}{\eps+\eps^{3/2}}\right| \ll \eps^{\theta -1} ,
\end{equation*}
thus the integers $Q_j^+=\big[ \frac{\cos
\omega_j}{\eps-\eps^{3/2}}\big] +1$ and $Q_j^-=\big[ \frac{\cos
\omega_{j+1}}{\eps+\eps^{3/2}} \big]$ satisfy
\begin{equation}\label{43}
0< Q_j^+-Q_j^- \ll \eps^{\theta-1}
\end{equation}
and
\begin{equation}\label{44}
\frac{1}{Q_j^+} \leq \frac{\eps -\eps^{3/2}}{\cos \omega_j} \leq
\frac{\eps+\eps^{3/2}}{\cos \omega_{j+1}} \leq \frac{1}{Q_j^-} \,
.
\end{equation}
Furthermore, it follows from \eqref{6.2} that there exists $\eps_0
=\eps_0 (\theta)>0$ such that for all $\eps <\eps_0$, all $j$, and
all $\frac{a}{q} \in [\tan \omega_j,\tan \omega_{j+1}]$, one has
\begin{equation*}
\frac{\eps -\eps^{3/2}}{\cos \omega_j} \leq \eps_\pm (q,a) \leq
\frac{\eps+\eps^{3/2}}{\cos \omega_{j+1}} \, .
\end{equation*}
This implies in conjunction with \eqref{44}, for all $\frac{a}{q}
\in [\tan \omega_j,\tan \omega_{j+1}]$, the inequalities
\begin{equation}\label{45}
\frac{1}{Q_j^+} \leq \eps_\pm (q,a) \leq \frac{1}{Q_j^-} \, .
\end{equation}
Since $Q_j^\pm =\frac{\cos \omega_j}{\eps}+O(\eps^{\theta-1})$,
one has
\begin{equation}\label{101}
(Q_j^\pm)^r=\frac{\cos^r \omega_j}{\eps^r}+O_r (\eps^{-r+\theta}).
\end{equation}

The first exit time increases when all the sizes of scatterers
decrease. Thus we infer from \eqref{45} the inequalities
\begin{equation}\label{46}
\int\limits_{I_j} l^r_{\frac{1}{Q_j^-}} (\omega)\, d\omega \leq
\int\limits_{I_j} \tilde{\tilde{\tau}}_\eps^r (\omega)\, d\omega
\leq \int\limits_{I_j} l^r_{\frac{1}{Q_j^+}} (\omega)\, d\omega .
\end{equation}
But by Theorem \ref{T1.2} and by \eqref{101} we may write
\begin{equation}\label{47}
\int\limits_{I_j} l^r_{\frac{1}{Q_j^\pm}} (\omega)\, d\omega= c_r
(Q_j^\pm)^r\int\limits_{I_j}\frac{dx}{\cos^r x}+ O_{r,\delta}
(\eps^{-r+\frac{1}{2}-2\alpha-\delta}+\eps^{-r+\theta+\alpha})
\end{equation}
with the better error term $\eps^{-\frac{3}{2}-\delta}$ for $r=2$.
Also using $\int_{I_j} \frac{dx}{\cos^r x} \ll_r \vert I_j \vert
\ll \eps^\theta$ we infer that the first integral in \eqref{47} is
expressible as
\begin{equation}\label{48}
\frac{c_r \cos^r \omega_j}{\eps^r}
\int\limits_{\omega_j}^{\omega_{j+1}} \frac{dx}{\cos^r x}
+O_{r,\delta}
(\eps^{-r+2\theta}+\eps^{-r+\frac{1}{2}-2\alpha-\delta}
+\eps^{-r+\theta+\alpha}),
\end{equation}
with the better error term
$\eps^{-2+\theta}+\eps^{-\frac{3}{2}-\delta}$ for $r=2$.

Summing up over $j$ we infer from \eqref{46}, \eqref{48} and
\eqref{42} that
\begin{equation}\label{49}
\begin{split}
\int\limits_I \tilde{\tau}^r_\eps (\omega)\, d\omega & =
\frac{c_r}{\eps^r} \sum\limits_{j=1}^N \cos^r \omega_j
\int\limits_{\omega_j}^{\omega_{j+1}} \frac{dx}{\cos^r x} \\
& \qquad +O_{r,\delta}
(\eps^{-r+\frac{1}{2}-2\alpha-\theta-\delta}
+\eps^{-r+\alpha}+\eps^{-r+\theta})
\end{split}
\end{equation}
with the better error term $\eps^{-\frac{3}{2}-\theta-\delta}
+\eps^{-2+\theta}$ for $r=2$.

Finally, we apply the mean value theorem and chose some $\xi_j \in
[\omega_j,\omega_{j+1}]$ to evaluate the sum
\begin{equation*}
\sum\limits_{j=1}^N \cos^r \omega_j
\int\limits_{\omega_j}^{\omega_{j+1}} \frac{dx}{\cos^r x}
\end{equation*}
as
\begin{equation*}
\begin{split}
\sum\limits_{j=1}^N \frac{(\omega_{j+1}-\omega_j)\cos^r
\omega_j}{\cos^r \xi_j} & =\sum\limits_{j=1}^N
(\omega_{j+1}-\omega_j)
\Big( 1+O_r (\omega_{j+1}-\omega_j)\Big) \\
& =\sum\limits_{j=1}^N (\omega_{j+1}-\omega_j) \Big(
1+O_r (\eps^\theta)\Big) \\
& =\vert I\vert +O_r (\eps^\theta) .
\end{split}
\end{equation*}
This implies Theorem \ref{T1.1} in conjunction with \eqref{49} by
taking $\theta=\alpha=\frac{1}{8}$ for $r\neq 2$ and
$\theta=\frac{1}{4}$ for $r=2$.
\end{proof}

\begin{proof}[Proof of Theorem \ref{T1.4}]
We proceed along the same line to estimate the moments of
$\tilde{R}$. Here we denote by $\tilde{\tilde{R}} (\omega)$ the
number of reflections in the side cushions in the case of vertical
scatterers (of variable size) $\{ q\} \times [a-\eps_-
(q,a),a+\eps_+(q,a)]$, $(q,a)\in \Z^{2\ast}$. It is seen as in the
proof of Theorem \ref{T1.1} that $\int_I \tilde{R}^r_\eps (\omega)
\, d\omega$ differs from $\int_I \tilde{\tilde{R}}^r_\eps
(\omega)\, d\omega$ by an error term of order $O_r (\eps^{2-r})$.
One can also show that
\begin{equation}\label{100}
\int\limits_{I_j} R^r_{\frac{1}{Q_j^-}} (\omega)\, d\omega \leq
\int\limits_{I_j} \tilde{\tilde{R}}_\eps^r (\omega)\, d\omega \leq
\int\limits_{I_j} R^r_{\frac{1}{Q_j^+}} (\omega) \, d\omega .
\end{equation}
Applying now Theorem \ref{T1.5} to the vertical scatterers
$V_{1/Q_j^\pm}$ on the intervals $I_j= [\omega_j,\omega_{j+1}]$ of
equal size $\vert I_j \vert =\frac{\vert I\vert}{N} \asymp
\eps^{\frac{1}{8}}$ with $\theta=\alpha=\frac{1}{8}$, and also
using \eqref{101}, we find that
\begin{equation*}
\int\limits_{I_j} R^r_{\frac{1}{Q_j^\pm}} (\omega)\, d\omega
=\frac{c_r \cos^r \omega_j}{\eps^r} \int\limits_{I_j} (1+\tan
x)^r\, dx + O_{r,\delta} (\eps^{-r+\frac{1}{4}-\delta}),
\end{equation*}
and thus
\begin{equation}\label{50}
\int\limits_I \tilde{R}^r_\eps (\omega) \, d\omega =
\frac{c_r}{\eps^r} \sum\limits_{j=1}^N \cos^r \omega_j
\int\limits_{\omega_j}^{\omega_{j+1}} (1+\tan x)^r \, dx+
O_{r,\delta} (\eps^{-r+\frac{1}{8}-\delta}) .
\end{equation}

By the mean value theorem we find $\xi_j,\eta_j \in I_j$ such that
\begin{equation}\label{51}
\sum\limits_{j=1}^N \cos^r \omega_j
\int\limits_{\omega_j}^{\omega_{j+1}} (1+\tan x)^r \,
dx=\sum\limits_{j=1}^N (\omega_{j+1}-\omega_j)\cos^r \omega_j
(1+\tan \xi_j)^r,
\end{equation}
and respectively
\begin{equation*}
\begin{split}
\int\limits_I & (\sin x+\cos x)^r \, dx =\sum\limits_{j=1}^N
\int\limits_{\omega_j}^{\omega_{j+1}} (\sin x+\cos x)^r \, dx  \\
& =\sum\limits_{j=1}^N (\omega_{j+1}-\omega_j)(\sin \eta_j+\cos
\eta_j)^r =\sum\limits_{j=1}^N (\omega_{j+1}-\omega_j) \cos^r
\eta_j (1+\tan \eta_j)^r .
\end{split}
\end{equation*}

From
\begin{equation*}
\cos^r \omega_j =\cos^r \eta_j +O_r (\omega_{j+1}-\omega_j)=
\cos^r \eta_j +O_r (\eps^{\frac{1}{8}})
\end{equation*}
and
\begin{equation*}
(1+\tan \xi_j)^r =(1+\tan \eta_j)^r +O_r (\vert \tan \xi_j -\tan
\eta_j \vert )= (1+\tan \eta_j)^r+O_r (\eps^{\frac{1}{8}})
\end{equation*}
we infer that the sum in \eqref{51} is equal to
\begin{equation*}
\begin{split}
\sum\limits_{j=1}^N (\omega_{j+1}-\omega_j) & \Big( \cos^r
\eta_j+O_r (\eps^{\frac{1}{8}})\Big) \Big(
(1+\tan \eta_j)^r +O_r (\eps^{\frac{1}{8}})\Big) \\
& =\int\limits_I (\sin x+\cos x)^r \, dx+O_r (\eps^{\frac{1}{8}}).
\end{split}
\end{equation*}
This can be combined with \eqref{50} to collect
\begin{equation*}
\int\limits_I \tilde{R}^r_\eps (\omega)\, d\omega =
\frac{c_r}{\eps^r} \int\limits_I (\sin x+\cos x)^r \, dx+
O_{r,\delta} (\eps^{-r+\frac{1}{8}-\delta}),
\end{equation*}
which concludes the proof of Theorem \ref{T1.4}.
\end{proof}

\end{document}